\documentclass[12pt]{article}
\usepackage{amssymb}
\usepackage{amsthm}
\usepackage{amsxtra}
\usepackage{amsfonts}
\usepackage{xr}

\usepackage{color}

\begin{document}

\def\sect{\section}

\newtheorem{thm}{Theorem}[section]
\newtheorem{cor}[thm]{Corollary}
\newtheorem{lem}[thm]{Lemma}
\newtheorem{prop}[thm]{Proposition}
\newtheorem{propconstr}[thm]{Proposition-Construction}
\newtheorem{fact}[thm]{Fact}

\theoremstyle{definition}
\newtheorem{para}[thm]{}
\newtheorem{ax}[thm]{Axiom}
\newtheorem{conj}[thm]{Conjecture}
\newtheorem{defn}[thm]{Definition}
\newtheorem{notation}[thm]{Notation}
\newtheorem{rem}[thm]{Remarks}
\newtheorem{remark}[thm]{Remark}
\newtheorem{question}[thm]{Question}
\newtheorem{example}[thm]{Example}
\newtheorem{problem}[thm]{Problem}
\newtheorem{excercise}[thm]{Exercise}
\newtheorem{ex}[thm]{Exercise}

\def\Bbb{\mathbb}
\def\cal{\mathcal}
\def\mL{{\mathcal L}}
\def\mC{{\mathcal C}}

\overfullrule=0pt

\def\si{\sigma}
\def\prf{\smallskip\noindent{\it        Proof}. }
\def\call{{\cal L}}
\def\nat{{\Bbb  N}}
\def\la{\langle}
\def\ra{\rangle}
\def\inv{^{-1}}
\def\ld{{\rm    ld}}
\def\trdeg{{\rm tr.deg}}
\def\dim{{\rm   dim}}
\def\th{{\rm    Th}}
\newcommand{\rest}[1]{{\lower       .25     em      \hbox{$\vert #1$}}}
\def\ch{{\rm    char}}
\def\zee{{\Bbb  Z}}
\def\conc{^\frown}
\newcommand\acl{{\rm acl}_\si}
\def\cls{{\rm cl}_\si}
\def\cals{{\cal S}}
\def\mult{{\rm  Mult}}
\def\calv{{\cal V}}
\def\aut{{\rm   Aut}}
\def\ffi{{\Bbb  F}}
\def\ffiti{\tilde{\Bbb          F}}
\def\degs{deg_\si}
\def\calx{{\cal X}}
\def\gal{{\cal G}al}
\def\cl{{\rm cl}}
\def\loc{{\rm locus}}
\def\calg{{\cal G}}
\def\calq{{\cal Q}}
\def\calr{{\cal R}}
\def\caly{{\cal Y}}
\def\aff{{\Bbb A}}
\def\cali{{\cal I}}
\def\calu{{\cal U}}
\def\epsilon{\varepsilon} 
\def\Uu{{\cal U}}
\def\rat{{\Bbb Q}}
\def\ga{{\Bbb G}_a}
\def\gm{{\Bbb G}_m}
\def\cee{{\mathbb C}}
\def\ree{{\mathbb R}}
\def\frob{{\rm Frob}}
\def\Frob{{\rm Frob}}
\def\fix{{\rm Fix}}
\def\Uu{{\cal U}}
\def\proj{{\Bbb P}}
\def\sym{{\rm Sym}}
 
\def\dcl{{\rm dcl}}
\def\calm{{\mathcal M}}

\font\helpp=cmsy5
\def\semdp
{\hbox{$\times\kern-.23em\lower-.1em\hbox{\helpp\char'152}$}\,}

\def\dnfo{\,\raise.2em\hbox{$\,\mathrel|\kern-.9em\lower.35em\hbox{$\smile$}
$}}
\def\dnf#1{\lower1em\hbox{$\buildrel\dnfo\over{\scriptstyle #1}$}}
\def\dfo{\;\raise.2em\hbox{$\mathrel|\kern-.9em\lower.35em\hbox{$\smile$}
\kern-.7em\hbox{\char'57}$}\;}
\def\df#1{\lower1em\hbox{$\buildrel\dfo\over{\scriptstyle #1}$}}        
\def\stab{{\rm Stab}}
\def\qfcb{\hbox{qf-Cb}}
\def\perf{^{\rm perf}}
\def\sipm{\si^{\pm 1}}

\def\vlabel{\label}

\title{On subgroups of semi-abelian varieties defined by difference equations}

\author{Zo\'e Chatzidakis\thanks{partially supported by PITN-2009-238381
    and by ANR-06-BLAN-0183, ANR-09-BLAN-0047, ANR-13-BS01-0006.}{\ } --  CNRS (UMR 8553) - Ecole
    Normale Sup\'erieure 
\and 
Ehud Hrushovski\thanks{The research leading
to these results has received funding from the European Research Council under the European
Unions Seventh Framework Programme (FP7/2007- 2013)/ERC Grant Agreement No. 291111, as well as the
Israel Science Foundation 1048/07.}{\ } -- The Hebrew University of Jerusalem}
\date{}
%\centerline{\today}
\maketitle

\newcommand{\Aut}{{\rm Aut}}
\newcommand{\End}{{\rm End}}
\newcommand{\calo}{{\cal O}}
\newcommand{\calp}{{\cal P}}
%\newcommand{\calm}{{\cal M}}
%\newcommand{\calq}{{\cal Q}}
% \font\helpp=cmsy5
% \def\semdp{\hbox{$\times\kern-.23em\lower-.1em\hbox
% {\helpp\char'152}$}\,}
% \def\calh{{\cal H}}
% \def\dd{dd}
% \def\abs#1{\vert #1 \vert}
% \def\ga{{\Bbb G}_a}
% \def\gm{{\Bbb G}_m}
\newcommand{\SU}{{\rm SU}}
\newcommand{\evSU}{{\rm evSU}}
\newcommand{\evsu}{{\rm evSU}}
\newcommand{\Sing}{{\rm Sing}}
\newcommand{\Fr}{{\rm Frob}}
\newcommand{\Supp}{{\rm Supp}}
\newcommand{\co}{{\rm co}}
\newcommand{\Hom}{{\rm Hom}}
\newcommand{\Ker}{{\rm Ker \,}}
\newcommand{\Fix}{{\rm Fix}}
\newcommand{\GL}{{\rm GL}}
\newcommand{\Zz}{{\mathbb Z}}
\newcommand{\Ss}{{\mathbb S}}

\section*
{\bf Introduction}

Consider the algebraic dynamics on an algebraic torus $T=\gm^n$ given by a matrix $M \in \GL_n(\Zz)$. 
 Assume no root of unity is   an eigenvalue of  $M$.  We show that any finite, equivariant map from another 
algebraic dynamics into $(T,M)$ arises from a group isogeny $\gm^n \to
\gm^n$ (see \ref{sketch-1} and more generally \ref{dyn}).    In other words,
the automorphism  $x \mapsto x^M$ of $K(x)=K(x_1,\ldots,x_n)$ does not extend to any finite field extension,
except those contained in $K(x_1^{1/m},\ldots,x_n^{1/m})$ for some
 $m\geq 1$.  
 A similar statement is shown for abelian varieties,
and in fact for semi-abelian varieties.

More generally, we study irreducible difference equations of the form  $n \si(x)=Mx$, with $M \in \End(A)$, $n \in \nat$; for instance the equation 
  $\si(x)^3=x^2$ on $\gm$.     We obtain a similar statement for the function field of such equations.

  Model-theoretically, this completes the description (\cite{[CH]}, \cite{[CHP]}, \cite{[H]})  of the induced structure on ACFA-definable subgroups of semi-abelian
  varieties.  Such subgroups (up to finite index) are defined by difference equations of the form 
  $n \si(x)=Mx$, with $M \in \End(A)$.  The induced structure is stable, except when the equation involves the points $A(F)$
of the fixed field or a twisted fixed field $\si^r(x)=x^{p^m}$.    Whereas the quantifier-free induced structure was understood previously -- it corresponds
to invariant subvarieties -- the full induced structure involves also finite covers, and stability was known only in characteristic zero.  
%Note that stability
%fails for the additive group in positive characteristic. 

We proceed to describe the result in terms of difference algebra.
By a difference field, we mean a field $K$ with a
distinguished automorphism $\si$. 
The theory ACFA of existentially closed difference
fields  was
extensively studied in \cite{[CH]} and \cite{[CHP]}. In these papers, a
characterisation of modular types was given, and it was shown that,
in characteristic $0$, all modular types are stable and stably
embedded. In characteristic $p>0$, we
however exhibited examples of modular subgroups of the additive group
$\ga$ which are not stable. The main result of this paper,
Theorem \ref{thm1}, is that all 
modular definable subgroups of a semi-abelian variety are stable and stably
embedded. This implies that definable subsets of  modular subgroups are Boolean combinations
of cosets of definable subgroups. 

Let $F$ be the transformal function field of a definable subgroup $B$ of
a semi-abelian variety $A$.  
Stability of $B$ is equivalent to the existence of few difference field
extensions $L$ of $F$, that are finite as field extensions.  In fact we
obtain a complete description of such extensions.  In characteristic
zero,
the main geometric tool is the ramification divisor of $L$ over
appropriate
varieties ($A$ or powers of $A$.)  With controlled exceptions, the
ramification divisor of a potential extension $L$ is invariant under the
dynamics, leading to reduction of the dimension of $B$.  For abelian
varieties, the ramification divisor still carries enough
information. If $B$ lives on a torus $\gm^n$ however, too many finite
extensions have the same ramification divisor, and a finer invariant
is needed.  We consider a certain invariant subspace of the Berkovich space, consisting of
valuations with center contained in the ramification divisor.  We define 
an invariant, in the value group, associated with a wildly ramified
extension $L$ of $F$.  Within this subspace, there may be no fixed
points but there are always recurrent points of the dynamics.  Such
points lead again to an eigenvector of the dynamics (acting on the value
group now) and to a reduction of the dimension of $B$.    {Such invariants of wildly ramified extensions are new to our
knowledge, and may be of interest elsewhere.  The use of recurrent points (rather than periodic points) also seems noteworthy.}

\smallskip
The paper is organised as follows. In section 1 we set up the notation and
recall some results on existentially closed difference fields from \cite{[C]}, \cite{[CH]} and
\cite{[CHP]}. Section 2 recalls the tools used to study definable subgroups of
algebraic groups, and describes the criterion for modularity of a
definable subgroup of a simple abelian variety or of $\gm$. 

Section 3 contains a host of technical lemmas, which will be used in the
proofs of Propositions \ref{ram3b} and \ref{prop1}. Section 4 contains the
proof of Theorem \ref{thm1}, and derives some consequences.

%% Proposition \ref{prop1} is due to the first author. The rest of the paper
%% was worked out conjointly.
{We would like to thank the referee for many useful comments.}

\section{Notation, preliminary definitions and results}

\para {\bf Notation and conventions}. \vlabel{notation} We work in the language
$\call=\{+,-,\cdot, 0,1,\si\}$ of difference fields. {\bf Unless 
   explicitly stated otherwise, all difference
fields are inversive}, i.e., the endomorphism $\si$ is surjective.
{\bf Throughout the paper,  we work inside a large
saturated model $(\Omega,\si)$ of ACFA}. If $K$ is a difference
subfield of $\Omega$, and $A\subset \Omega$, $K(A)_\si$ denotes the
difference field generated by $A$ over $K$,
 and $\acl(A)$ the smallest algebraically closed field containing $A$ and
closed under $\si$ and $\si\inv$. If $A$ is a subfield
of $\Omega$, then $A^{alg}$ denotes the (field-theoretic) algebraic
closure  of $A$ and $A^s$ its separable closure. We let $\Fr$ denote the identity of $\Omega$ if the
characteristic of $\Omega$ is $0$, and the automorphism $x\mapsto x^p$
if the characteristic is $p>0$. If $char(\Omega)=p>0$ and $q$ is a
power of $p$, we also denote by $\Frob_q$ the automorphism $x\mapsto
x^q$ of $\Omega$. If $L$ is a finite algebraic extension of $K$, then
$[L:K]$ denotes its degree, $[L:K]_s=[L\cap K^s:K]$ its {\em separable
  degree}, and $[L:K]_i=[L:K]/[L:K]_s$ its {\em inseparable
  degree}. \\[0.1in]
If $n$ is a positive integer, we denote by $\call[n]$ the language
$\{+,-,\cdot, 0,1,\si^n\}$, viewed as a sublanguage of $\call$, and
by $\Omega[n]$ the difference field $(\Omega,\si^n)$. Then
$\Omega[n]$ is a model of ACFA (\cite{[CH]}, Corollary (1.12)(1)), and is
saturated. If $E$
is a difference subfield of $\Omega$, and $a$ a tuple of $\Omega$,
then $tp(a/E)[n]$ denotes the type of $a$ over $E$ in the structure
$\Omega[n]$.\\[0.1in]
Recall that the ring of difference polynomials over $K$ in $\bar
X=(X_1,\ldots,X_n)$, denoted $K[\bar X]_\si$, is simply the polynomial ring
$K[\si^j(X_i)\mid i=1,\ldots,n,\, j\in\nat]$. The
{\it $\si$-topology on $K^n$} is the topology with basic closed
sets the {\it $\si$-closed sets} $\{\bar a\in K^n\mid f(\bar a)=0\}$,
where $f(\bar X)$
is a tuple of difference polynomials over $K$. This topology is
Noetherian (see \cite{[C]}, 3.V).
When working inside $\Omega[n]$ we will speak of the
$\si^n$-topology. \\[0.1in]
If $V$ is an algebraic set, and $\tau$ an automorphism, we will denote by $\tau(V)$ or by $V^\tau$ the algebraic set whose defining equations are obtained by applying $\tau$ to the coefficients of the equations defining $V$. In particular, $V^\tau(\Omega)=\tau(V(\Omega))$. 

\para\vlabel{ACFA-1} {\bf Basic notions and results on ACFA}. Model-theoretic algebraic closure coincides
with  $\acl$
 ((1.7) in \cite{[CH]}), and (model-theoretic)
independence of algebraically closed sets $A$ and $B$ over a common
algebraically closed subset $C$ corresponds to linear independence of
the fields $A$ and $B$ over $C$, and any completion of ACFA is supersimple
(see (1.9) in \cite{[CH]} and use \cite{[KiP]}). We also know that any completion of ACFA
eliminates imaginaries ((1.12) in \cite{[CH]}).

\smallskip\noindent
Let
$E$ be a difference subfield of $\Omega$, and $a$ a tuple of elements
of $\Omega$. Then the quantifier-free type of $a$ over $E$, denoted
$qftp(a/E)$, is the set of quantifier-free $\call$-formulas with parameters
in  $E$ which are satisfied by $a$. It therefore describes the
isomorphism type of the difference field $E(a)_\si$ over $E$. Similarly,
$qftp(a/E)[n]$ denotes the set of quantifier-free $\call[n]$-formulas with
parameters in $E$ which are satisfied by $a$.

\smallskip\noindent
The SU-rank is defined as usual. Let us mention that $\SU(a/A)$ is
finite if and only if \\${\rm tr.deg}(\acl(A,a)/\acl(A))$ is finite, if and
only if all elements of the tuple $a$ satisfy some
non-trivial difference equation over $\acl(A)$. We denote by
$\SU(a/A)[n]$ the SU-rank in the reduct $\Omega[n]$ (thus it equals
$\SU(a/acl_{\si^n}(A))[n]$).\\[0.1in]
Finally, assume that $A$ is a difference subfield of
$\Omega[n]$ for some $n$, and let $a\in \Omega$. We define the
eventual SU-rank of $a$ over $A$, denoted $\evSU(a/A)$, as
$lim_{m\rightarrow \infty}\SU(a/A)[m!]$ (see (1.13) in \cite{[CHP]}). It is
well-defined, and only depends on $qftp(a/A)$. Note that if $A=\acl(A)$ and
$\si(a)\in A(a)^{alg}$, then
$\SU(a/A)[m]\leq \SU(a/A)[mn]$ for all $m,n\neq 0$. This implies that
there is some $m>0$ such that for all $n>0$,
$\evSU(a/A)=\SU(a/A)[mn]$.

\para\vlabel{bab}{\bf Completions of  quantifier-free types}. Let
$E$ be a difference subfield of $\Omega$, and $a$ a tuple of elements
of $\Omega$. Then $tp(a/E)=tp(b/E)$ if and only if there is an
isomorphism $\acl(E,a)\to \acl(E,b)$, which leaves $E$ fixed and
sends $a$ to $b$ ((1.15) in \cite{[CH]}). In particular, if $E=\acl(E)$, then
ACFA$\,\cup qftp(E)$ is complete, and the completions of ACFA are
obtained by describing the isomorphism type of the algebraic closure
of the prime field.\\[0.05in]
A {\it finite $\si$-stable extension} of a difference subfield $K$
of $\Omega$ is a finite (algebraic) extension $L$ of $K$ such that
$\si(L)=L$. If $L$ is a finite separable $\si$-stable extension of $K$, then so is its Galois
closure $M$ over $K$. Furthermore, whether or not $\si(M)=M$ does not
depend on the extension of $\si$ to $M$, but is completely determined by
the isomorphism type of $K$: if $M$ is finite separable over $K$, let $\alpha\in
M$ be such that $M=K(\alpha)$, and let $P(T)\in K[T]$ be the minimal
monic polynomial of $\alpha$ over $K$, and $P^\si(T)$ the polynomial
obtained by applying $\si$ to the coefficients of $P$. Then $\si(M)=M$
is equivalent to the following statement: the field $K[T]/(P(T))$  contains
a root  of $P^\si(T)$.  Observe also that $\si$ extends uniquely to
the perfect hull of $K$.   

\begin{thm}\vlabel{bab2} (Babbitt, \cite{[C]} Theorem 7.VIII). Let $E$
be a difference subfield of $\Omega$, and $a$, $b$ tuples which have
the same quantifier-free type over $E$. The following conditions are
equivalent:
\begin{enumerate}
\item{$tp(a/E)=tp(b/E)$.}
\item{Given any finite $\si$-stable Galois extension $L$    of
$E(a)_\si$,   there is an $E$-embedding $L\to
\Omega$ which sends $a$ to $b$.}
\end{enumerate}
\end{thm}
\noindent
In particular, if $E(a)_\si$ has no non-trivial finite separable
$\si$-stable extension, then $qftp(a/E)$ is complete.

\smallskip\noindent Let $L$ be a finite separable $\si$-stable Galois extension of
$E(a)_\si$. Then $\si$ induces an automorphism of
$G=\gal(L/E(a)_\si)$ given by $\rho\mapsto \si\inv\rho\si$; hence, for
some $\ell$ we will have $\si^{-\ell}\rho\si^\ell=\rho$ for all $\rho\in
G$. Difference fields between $E(a)_\si$ and $L$ correspond to subgroups
$H$ of $G$ which are stable under the action of $\si$. \\[0.1in]
While it may
happen that $E(a)_\si$ has some non-trivial
finite separable $\si$-stable extension, and yet $qftp(a/E)$ be complete, this
does not hold if one wants to consider all the reducts $\Omega[n]$.
Namely ((2.9)(5) in \cite{[CH]}):

\smallskip\noindent
\begin{lem}\vlabel{bab3} Let $F$
be a difference subfield of $\Omega$. The following
conditions are equivalent:
\begin{enumerate}
\item{ACFA $\cup \,qftp(F)[n]$ is complete, for all $n>0$.}
\item{$F$ has no finite separable $\si$-stable extension.}
\end{enumerate}
\end{lem}

\para\vlabel{mod1}{\bf Modularity}. Let $S\subset \Omega^n$ be stable under any
$\acl(E)$-automorphism of $\Omega$. We say that $S$ is {\em modular} if
whenever $a$ is a tuple of elements of $S$ and
$B\subset \Omega$, then $a$ and $B$ are independent over
$\acl(E,a)\cap \acl(E,B)$. A type over $E$ is {\it modular} if the set
of its realisations is modular. \\[0.1in]
This notion of modularity is also called {\it one-basedness}. Here we
use the fact that our theory is supersimple and eliminates imaginaries.
The main result of \cite{[CHP]} (see (7.5)) shows that $S$ is modular if and only
if, for
every $F=\acl(F)$ containing $E$ and $a\in S$, $tp(a/F)$ is orthogonal
to all fixed fields, i.e., to all formulas of the form
$\si^n(x)=\Frob^m(x)$ where $n\geq 1$, $m\in\zee$. It follows that
  $tp(a/E)$ is
modular if and
only $tp(a/E)[\ell]$ is modular for every $\ell\geq 1$, and that a modular
type has finite $\SU$-rank. 

\begin{thm}\vlabel{f1} {\rm (\cite{[CHP]} (7.6))}. Let $G$ be an algebraic group
and
$B$ a definable modular subgroup of $G(\Omega)$ of finite SU-rank. If $D$ is a
quantifier-free definable subset of $B$, then $D$ is a finite Boolean
combination of translates of quantifier-free definable subgroups of
$B$. If $B$ is defined over $E$, then every quantifier-free definable
subgroup of $B$ is definable over $\acl(E)$.
\end{thm}

\para\vlabel{sse}{\bf Stability and stable embeddability}.
A subset $S$ of $\Omega^m$ stable
under $E$-automorphisms of $\Omega$ is {\em stably embedded} if whenever
$D\subset \Omega^{nm}$ is definable, then $D\cap S^n$ is definable with
parameters from $S$ (see the
Appendix of
\cite{[CH]} for properties). A type is {\it stably embedded} if the set of
its realisations is stably embedded. \\[0.1in]
Let $E=\acl(E)$ and $a$ a 
tuple of elements of $\Omega$. If $tp(a/E)$ is stationary (i.e.,
has a unique non-forking extension to any set $F\supset E$), then
$tp(a/E)$ is definable (this is standard, using compactness).  Hence, if
all extensions of
$tp(a/E)$ to algebraically closed sets are stationary,
then $tp(a/E)$ is {\em stable and stably embedded}
(see Lemma 2 and its Definition-Remark in the Appendix of
\cite{[CH]}). A definable set $D$ is {\em stable and stably embedded} if all
types realised in $D$ are stable and stably embedded. Equivalently, if $D$
with the structure induced by $\Omega$ is stable and stably embedded. \\[0.1in]
Using the
characterisation of modular types and
Lemmas 2, 3 of the Appendix of \cite{[CH]}, one easily deduces:

\smallskip\noindent
\begin{prop}\vlabel{sse2} Let
$$1\longrightarrow B_1 \longrightarrow B_2 \longrightarrow B_3
\longrightarrow 1$$ be  an exact
sequence of groups definable in a model of ACFA. Then $B_2$ is
modular [resp. stable and stably embedded] if and only if $B_1$ and
$B_3$ are modular [resp. stable and stably embedded].
\end{prop}

\section{\bf Definable subgroups of algebraic groups}
In this chapter we introduce tools used to study definable groups, and
give a brief sketch of the description of modular subgroups of
abelian varieties. The proof in \cite{[H]} was given in characteristic  $0$, and we
indicate what changes need to be made in positive characteristic. Please
see chapter 4 in \cite{[H]} for more details.

\para\vlabel{gp1}{\bf Setting and notations.} In what follows we will
use $p$ to denote $char(\Omega)$ if it is positive, and $1$ if it is
$0$. Thus Frob will denote the map $x\mapsto x^p$ (see \ref{notation}). Let $G$ be a connected algebraic group, $H$ a definable subgroup of
$G(\Omega)$, everything defined over $E=\acl(E)\subset \Omega$.\\[0.1in]
For  $m\in\nat$, define $G_{(m)}=G\times \si(G)\times \cdots\times
\si^m(G)$ and $p_m:G\to G_{(m)}$ by
$g\mapsto (g,\si(g),\ldots,\si^m(g))$. Let $H_{(m)}$ be
the Zariski closure of $p_m(H)$, and let $\tilde H_{(m)}=\{g\in G\mid
p_m(g)\in H_{(m)}\}$. The intersection $\tilde H$ of the $\tilde
H_{(m)}$ equals the $\si$-closure of $H$  (in $G(\Omega)$), and equals some
$\tilde H_{(m)}$, because every descending sequence of $\si$-closed sets
stabilises. Then $[\tilde H:H]<\infty$: one  easily shows that the
generics of $H$ (in the model-theoretic sense) are those $h\in H$
such that $p_m(h)$ is a generic of the
algebraic group $H_{(m)}$ for every $m>0$. Hence a generic of
$H$ is a generic of $\tilde H$.

Because the
$\si$-topology is Noetherian, every $\si$-closed set $S$ can be
expressed as an irredundant union of $\si$-closed irreducible sets,
which are called the irreducible components of $S$. Then the
irreducible component  of $\tilde H$ containing the identity of $G$
is a subgroup of $\tilde H$; it is called the {\em connected component
of} $\tilde H$ and is denoted by $\tilde H^0$. Then $[\tilde H:\tilde
H^0]<\infty$. We let $H^0=H\cap \tilde H^0$, and call it the {\em connected
component} of $H$. Then also $[H:H^0]<\infty$.

\para\vlabel{cmin}{\bf c-minimal subgroups}. Let $B$ be a definable subgroup of some
algebraic group. Then $B$ is
{\em c-minimal} iff every definable subgroup of $B$ is either finite
or of finite index in $B$.  
Note that c-minimality is preserved by
definable homomorphisms with finite kernel; hence if $B$ is
c-minimal, then there is a definable homomorphism $g:B^0\to G(\Omega)$, where
 $\Ker(g)$ is finite central, and
$G$ is an  algebraic group which is simple  and is the Zariski closure of
$g(B^0)$. Indeed, one simply looks at the simple 
quotients  of the Zariski closure of $B^0$, and
at the images of $B^0$ in these quotients: one of these images must be
isogenous to $B^0$. \\
Note also that  if $SU(B)=1$ then $B$ is
c-minimal. If $B$ is modular, then the converse is true as well:
c-minimality of $B$ implies $SU(B)=1$: this is because the stabilizer of
any type realised in $B$ has the same SU-rank as the type.

\para \vlabel{end-0}{\bf Definable endomorphisms and subgroups of abelian varieties}.
Let $A$ be an abelian variety, defined over $E=\acl(E)$. By standard
results on abelian varieties, $A$ is isogenous to a finite direct sum
of simple abelian varieties defined over $E$, say to
$\bigoplus_{i=1}^n A_i$. Renumbering, we may assume that for $i<
j\leq r$, for all $\ell\in \zee$, $A_i$ is not isogenous to $\si^\ell(A_j)$,
and that for any $j>r$ there are $i\leq r$ and $\ell\in\zee$ such that
$A_i$ and $\si^\ell(A_j)$ are isogenous. For $i\leq r$ let $m(i)$ be the
number of indices $j$ such that for some $\ell$, $A_i$ and $\si^\ell(A_j)$
are isogenous.

\smallskip
Let us denote by $\End(A)$ the ring of (algebraic) endomorphisms of $A$,
by $\Hom(A,A_i)$ the group of algebraic homomorphisms $A\to A_i$, and
by $\End_\si(A)$ the ring of definable endomorphisms of $A(\Omega)$,
$\Hom_\si(A,A_i)$ the group of definable homomorphisms $A\to A_i$.
Hrushovski gives a good description of $\rat\otimes \End_\si(A)$ in
\cite{[H]}, and we refer to this paper for the results quoted below.  First,
note that 
$\rat\otimes \End_\si(A)\simeq \prod_{i=1}^r M_{m(i)}(\rat\otimes
\End_\si(A_i))$, so that it suffices to describe the rings
$\End_\si(A_i)$. We will therefore restrict our attention to simple
abelian varieties. \\[0.05in]
Recall that two definable subgroups $B$ and $C$ of a group are {\it
commensurable} if $B\cap C$ has finite index in both $B$ and $C$. % (One
% can extend this definition to $\infty$-definable subgroups by
% requiring that $B\cap C$ be of bounded index in $B$ and in $C$).
Hrushovski shows that  a definable subgroup of $A(\Omega)$ is
 commensurable to a definable subgroup $\bigcap_{j=1}^n\Ker(F_j)$,
where $F_j\in \Hom_\si(A,A_{j})$.

\smallskip\noindent
The study in \cite{[H]} is made for difference fields of characteristic $0$.
However, the proofs generalise to positive characteristic without any
trouble. Note that in positive characteristic, the Frobenius map may
define an endomorphism of the variety.

\medskip\noindent
\begin{thm}\vlabel{end-1}Let $A$ be a simple abelian
variety
defined over $E=\acl(E)$, let  $B$ be an infinite definable subgroup of
$A(\Omega)$.
\begin{enumerate}
\item If there is no integer $n>0$ such that $A$ and
$\si^n(A)$ are isogenous, then $\End_\si(A)=\End(A)$, and $B=A(\Omega)$.
\item Assume that there is an integer $n>0$ such that $A$
and $\si^n(A)$ are isogenous. We fix such an $n$,  smallest
possible, and choose an isogeny $h:A\to \si^n(A)$ of
minimal degree $m$. If $\si^n(A)=A$, then we choose $h$ to be the
identity. Let   $h':\si^n(A)\to A$ 
be such that $h'h=[m]$ (multiplication by $m$ in $A$; such an $h'$
exists by standard results on abelian varieties); then $hh'=[m]$.
Define  $\tau=\si\inv h$ and $\tau'=h'\si$.

Then $\rat \otimes \End_\si(A)\simeq \rat\otimes \End(A)[\tau,\tau']$, and $B$ is
commensurable to $\Ker(f)$ for some $f\in \End(A)[\tau,\tau']$.
Furthermore, $\rat \otimes \End_\si(A)$ is an Ore domain and if
$C\subseteq B$ is definable, then $C$ is commensurable to some $\Ker(g)$ with $g$
dividing $f$ (i.e, $hg=f$ for some $h\in \rat \otimes \End_\si(A)$). It
follows that $B$ is c-minimal if and only if $f$ is left-irreducible. 
\end{enumerate}
\end{thm}

\prf In characteristic $0$, this is given by Proposition 4.1.1 of
\cite{[H]}. The proof goes through verbatim.

\begin{thm}\vlabel{end-2}Let $A$ be a simple abelian
variety
defined over $E=\acl(E)$, and let $B=\Ker(f)$ be a {\bf proper}
definable subgroup of
$A(\Omega)$,  which is c-minimal
and connected. Assume that $B$ is  not
modular.

\begin{enumerate}
\item{Then $A$ is isomorphic to an abelian
variety $A'$ defined over $\Fix(\rho)$, where $\rho=\Frob ^m\si^n$ for
some $m\in\zee$, $n>0$.}

\item{Assume that $A$ is defined over $\Fix(\theta)$, where
$\theta=\Frob ^{s}\si^{t}$ for
some $s\in\zee$, $t>0$.  If $A$ is not isomorphic to any
variety defined over the algebraic closure of the prime field, then
$f$ divides $\theta^\ell -1$ for some $\ell$. If $A$ is isomorphic to a
variety $A'$ defined over the algebraic closure of the prime field, by
an isomorphism $F$, then
for some $m\in\zee$, $n>0$ and $\rho=\Frob ^m\si^n$, $F(B)\subset A'(\Fix(\rho))$. }
\end{enumerate}
\end{thm}

\prf First observe that the  c-minimality of $B$ implies that if $q$ is any
non-algebraic type realised in $B$, then the realisations of $q$
generate a subgroup of finite index in $B$. \\[0.05in]
(1) In characteristic $0$, this result is Proposition 4.1.2 of  \cite{[H]}. 
The proof  generalises easily to the positive characteristic case, using the
results of \cite{[CHP]}. Here is a sketch of the main steps: if $B$ is not
modular, then by c-minimality,  the type of a generic of $B$ is
non-orthogonal to one of the fixed fields, say $k=\Fix(\rho_0)$, where
$\rho_0=\si^n\Frob^m$, $n\geq 1$, and $(n,m)=1$ if $m\neq 0$. Thus,
modulo a finite kernel, it is qf-internal to $k$, i.e., there is a finite
subgroup $C$ of $B$, and a definable  map $g_0$ from some definable
set $S\subset k^\ell$ onto $B/C$ (in fact $g_0$ is given piecewise by
difference rational functions).

Elimination of imaginaries in ACFA tells us that $B/C$ is then definably
isomorphic (via some $g_1$) with a  group $H_0$ living in some cartesian power of
$k$. On the other hand, every subset of a cartesian power of $k$ which
is definable in $\Omega$, is already definable in the difference field $k$,
using maybe extra parameters from $k$ (see (7.1)(5) in \cite{[CHP]}). Note that $k$ has
SU-rank $1$ ((7.1)(1) in \cite{[CHP]}). An argument similar to the one given
in \cite{[HP2]} or in \cite{[KP]} then gives us a definable (in $k$) map $g_2:H_0'\to
H_1(k)$, where $\Ker(g_2)$ is finite, $H_0'$ is a subgroup of finite
index of $H_0$,  and $H_1$ is an algebraic group
defined over $k$. Then $H'_0\supseteq [N]H_0$ for some $N$, and
so $g_1([N]B/C)\subseteq H'_0$. Hence, replacing
$g_1$ by $g_1\circ[N]$, 
 we may assume that $H_0'=H_0$. (Recall that the $N$-torsion of all
 these groups is finite). 

Composing these maps, we therefore obtain a
definable group homomorphism $h:B\to H_1(k)$, with $\Ker(h)$ finite. We
may assume 
that the Zariski closure of $h(B)$ is all of $H_1$.

Since $B$ is $c$-minimal, so is $h(B)$, and this implies that if
$\pi:H_1\to A'$
is a projection of $H_1$ onto a simple quotient $A'$ of $H_1$,
then $\pi h(B)$
is Zariski dense in $A'$, and therefore $\Ker(\pi)\cap h(B)$ is finite.
Replacing $\rho_0$ by some power $\rho=\rho_0^r$ (and  $k$ by its
algebraic extension of degree $r$), we may assume that $A'$
and $\pi$ are defined over $k$.
Hence we get a definable map $B\to A'(k)$, with finite kernel. This
implies that $\Hom(A,\si^\ell(A'))\neq (0)$ for some $\ell$, and
applying some power of $\si$ to $A'$, we may assume that
$\Hom(A,A')\neq 0$. Thus $A'$ is a simple abelian variety, defined
over $k$, and isogenous to $A$. This implies that $A$ is isomorphic to
an abelian variety $A''$ defined over some finite extension of $k$.
Replacing $\rho$ by its appropriate power, we have shown (1).\\[0.05in]
(2) By the proof of (1), we have a definable map $\varphi:B\to A'(k)$,
with finite 
kernel $D$, where $\rho=\Frob ^{m}\si^{n}$ and
$k=\Fix(\rho)$. Moreover $A'$ is isogenous to
$A$. If $mt=sn$, then we may assume that $\theta=\rho$ and $A'=A$. If
$mt\neq sn$, then $A$ is isomorphic to a variety $A''$ defined over
$\Fix(\theta)^{alg}\cap \Fix(\rho)^{alg}=\ffi_p^{alg}$, and we are in the
second case. 
The result will follow from the following claim:

\smallskip\noindent
{\bf Claim}. Let $A$ be an abelian variety defined over $k=\Fix(\rho)$, let
$B$ be a definable subgroup of $A(\Omega)$ and
$\varphi:B\to A(k)$  a definable homomorphism with finite kernel
$D$. Then $B\subseteq A(\Fix(\rho^\ell))$ for some $\ell$.

\prf The graph of $\varphi$ is
a definable subgroup of $A^2$; as there are only
countably many of those (see 4.1.10 in \cite{[H]}), it must be defined over $k^{alg}$. 
Hence, $\varphi$ is definable over $k^{alg}$, and, replacing
$\rho$ by an appropriate power we
may assume that $\varphi$ is defined over $\Fix(\rho)$. Let $k_0\prec k$ be
such that everything is defined over $k_0$. Since $D$ is finite,
if $b\in B$ then $b\in \acl(k_0(\varphi(b)))=k_0(\varphi(b))_\si^{alg}$.
By compactness, there is $\ell$ such that
$[k_0(\varphi(b))_\si(b):k_0(\varphi(b))_\si]\leq \ell$ for every $b\in B$
and this implies that $b\in \Fix(\rho^{\ell!})$.

 \para \vlabel{end-3}{\bf Definable subgroups of tori}. Similar
results hold for tori, i.e., algebraic groups isomorphic to $\gm^n$
for some $n$. Recall that algebraic subgroups of $\gm^n$ 
 are defined by equations of the form
$\prod_{i=1}^n x_i^{m_i}=1$ for some integers $m_1,\ldots,m_n$.  These
subgroups are connected if the $m_i$'s have no common divisor. % The
% techniques of proof of Theorem \ref{end-1} then generalise to give
% $\rat\otimes\End_\si(\gm)\simeq \rat[\si,\si\inv]$, and
% $\rat\otimes \End_\si(\gm^n)\simeq M_n(\rat[\si,\si\inv])$, but we will
% not need this result.
Using the
description of definable subgroups of algebraic groups given in \ref{gp1}, 
one shows easily that a definable subgroup $B$ of $\gm$ is commensurable to
$\Ker(f(\si))$ for some $f(T)\in\zee[T]$, and that $B$ is c-minimal if and
only if $f$ is irreducible. Here, if $f=\sum m_iT^i$ and $a\in \gm$, we
define $f(\si)(a)=\prod \si^i(a)^{m_i}$. \\[0.05in] 
The proof of Theorem \ref{end-2} generalises to this setting to show  that modular subgroups of $\gm$ are those which
are commensurable to some
 $\Ker(f(\si))$, where $f(T)$ is relatively prime (in $\rat[T]$) to all elements of the
form $T^n-p^m$, with $n>0$ and $m\in\zee$. Indeed, let $B=\Ker(f(\si))$,
where $f\in\zee[T]$ and assume that  $B$ is not
modular.  Then for some $k=\fix(\si^n\Frob^m)$, there is a definable homomorphism
$h:B\to \gm(k)$ with infinite image. Hence $\Ker(h)$ is commensurable
with some $\Ker(g(\si))$ with $g(T)\in \zee[T]$, so that $g$ divides $f$ in
$\rat[T]$. The homomorphism $h$ then induces a homorphism $\pi:B/B\cap
\Ker(g(\si))\to \gm(k)$ with finite kernel. The proof of the claim in
\ref{end-2} generalises to show that $B/B\cap \Ker(g(\si))\subseteq
\gm(\fix(\si^{n\ell}\Frob^{m\ell}))$ for some $\ell$, so that $B$ is
commensurable to $\Ker((\si^{n\ell}-p^{m\ell})g(\si))$.

\begin{lem}\vlabel{end-4} Let $f \in \zee[T]$,
  $f(T)=\sum_{i=0}^n b_i T^i$, with $b_0b_n\neq 0$, and assume that the
  $b_i$'s are relatively prime. Let $a$ be a
  generic of $\Ker(f(\si))$ over some field $E=\acl(E)$\footnote{I.e., $a$ satisfies $\prod_{i=0}^n
    \si^i(a^{b_i})=1$, and $\trdeg(E(a)_\si/E)=n$.}, let $m\geq n$, and let $g \in \zee[T]$ be
  such that $g(\si)(a)=1$. 
\begin{enumerate}
\item
Then $g$ belongs to the ideal generated by $f$ in $\zee[T]$.
\item Write $g=\sum_{i=0}^m c_i T^i$. Then $b_n$ divides $c_m$. 
\item Let $r>1$, and let $h(T)\in\rat[T]$ be monic and of least
  degree such that, writing  $N$ for  the least positive integer such
  that $Nh(T)\in\zee[T]$, 
  then $Nh(\si^r)(a)=1$. Then every prime divisor of $N$ divides
  $b_n$. 
% \item Let $A\in {\rm GL}_n(\rat)\cap M_n(\zee)$, and let $b$ be a
%   generic solution of $\si(b)=Ab$ in $\gm^n(\calu)$. Assume that 
\item  ($char(\Omega)=p>0$) Let $m>0$ and consider $\tau=\si\frob^m$. Then $a$ is
  in 
  $\Ker s(\tau)$, where $s(T)= p^\ell\sum_{i=0}^n p^{-mi}b_i T^i  $,
  with $\ell =-\inf_i \{v_p(b_i)-mi\}$, $v_p$ being the $p$-adic valuation
  on $\rat$. Moreover, $s$ is irreducible.
\end{enumerate}

\end{lem}

\prf (1) The set of $h\in\rat[T]$ such that $h(\si)(a)=1$ is an
ideal. Since $\rat[T]$ is a Euclidian domain,  there are $q,r\in
\rat[T]$ such that $g=fq+r$ and $\deg(r)<\deg(f)=n$. If $r\neq 0$, 
then 
$g(\si)(a)=1$ implies $r(\si)(a)=1$, which implies
$\trdeg(E(a)_\si/E)\leq \deg(r)<\deg(f)$, a contradiction since we
assume $a$ to be a generic of $\Ker(f(\si))$. Hence $r=0$. Let $\ell$ be a
prime number, and define a valuation $v$ on
$\rat[T]$ extending the $\ell$-adic valuation $v_\ell$ by setting
$v(\sum a_iT^i)=\min_i\{v_\ell(a_i)\}$. Our assumption on $f$ implies
that $v_\ell(f)=0$, and our assumption on $g$ that $v(g)\geq 0$. It
follows that $v(q)\geq 0$. This being true for any prime $\ell$ implies that
$q\in\zee[T]$. \\[0.05in]
(2) Follows from (1).\\[0.05in]
%Let $\ell$ be a prime divisor of $b_n$. In the notation of (1), we
%need to show that $v_\ell(c_m)\geq v_\ell(b_n)$. But if
%$q=\sum_{i=0}^{m-n}d_i\si^i$, then $c_m=d_{m-n}b_n$. Hence, if
%$v_\ell(c_m)<v_\ell(b_n)$, then $v_\ell(d_{m-n})<0$, which contradicts
%(1). \\[0.05in]
(3) We now work in $\rat[T]$. Let $\alpha_1,\ldots,\alpha_n$ be the
roots of $f(T)$ in $\rat^{alg}$. Then the roots of $h(T)$  in
$\rat^{alg}$ are among $\alpha_1^r,
\ldots,\alpha_n^r$. % , so that 
% $$f(T)=b_n\prod_{i=1}^n (T-\alpha_i), \qquad
% h(T)=\prod_{i=1}^n(T-\alpha_i^r).$$
Let $\ell$ be a prime number not
dividing $b_n$. This means that all $\alpha_i$'s are integral algebraic
over $\zee_\ell$; hence so are all $\alpha_i^r$, and this implies that
the  coefficients of $h$ are in $\zee_\ell$, in other
words, that $\ell$ does not divide $N$. \\[0.05in]
(4) A simple calculation shows that $s(T)\in\zee[T]$ and  has relatively
prime coefficients, and that $a\in\Ker s(\tau)$. The roots of $s(T)$ are
$p^m\alpha_1,\ldots,p^m\alpha_n$, from which one deduces  the irreducibility of
$s(T)$. 

\section{Technical lemmas}

In this chapter we collect some technical lemmas which will be used in
the proof of the main results.

\para{\bf Valuations - basic results and notation}. We refer to the book
of Engler and Prestel \cite{EP} for all  notions and results. First a
definition: if  $E$ is a subring of the field $K$, then we say that
a valuation $v$ on $K$ is  an {\em $E$-valuation} if $v$ is trivial on
$E$. Given
a valuation 
$v$ on a field $K$, we denote by $\calo_v$ or $\calo_K$ its valuation
ring, by $\calm_v$ or $\calm_K$ its maximal ideal, by
$k_K$ its residue field $\calo_K/\calm_K$, and by $\Gamma(K)$
its value group. If $a\in\calo_v$, then we denote by $\bar a$ its
residue, i.e., the image of $a$ in $\calo_v/\calm_v$. The number $p$ will denote the {\em residual
  characteristic} of $v$ (i.e., the
characteristic of $k_K$) if it is positive, and $1$ if it is $0$.  Let
$L$ be a finite algebraic extension of $K$, and $w$ an extension of $v$
to $L$, $\Gamma(L)$ and $k_L$ the corresponding value group and residue
field. We define $e(w/v)=e(L/K)=[\Gamma(L):\Gamma(K)]$ (the {\em reduced ramification index} of
$w$ in $L$), $r(w/v)=r(L/K)=e(L/K)[k_L:k_K]_i$ (the {\em ramification index}) and
$f(w/v)=f(L/K)=[k_L:k_K]$ (the {\em residual degree} or {\em inertia degree} of $w$ in $L$). We say
that $w$ {\em ramifies over $K$}  if $r(L/K)>1$, and that $v$ {\em ramifies in} $L$ if
some extension of $v$ to $L$ ramifies over $K$.  

If
$L$ is a finite normal extension of $K$,  then all extensions of $v$ to
$L$ are conjugate, so that the numbers $e(L/K)$, $r(L/K)$ and $f(L/K)$ will not
depend on the choice of $w$. We denote by $g(L/K)$ the number of
extensions of $v$ to $L$. The numbers $e(L/K), f(L/K)$, $r(L/K)$ and
$g(L/K)$ are 
multiplicative in towers.

The {\em defect} of a finite algebraic extension $L$ of $K$ is the
number $d=d(L/K)$ such that $[L:K]=d\sum_{w}e(w/v)f(w/v)$ where $w$
 ranges over all extensions of $v$ to $L$. This number is a power of $p$
 (Theorem 3.3.3 in \cite{EP}). 
We say that the valued field $K$ is {\em defectless} if for every finite
 extension $L$ of $K$, we have $d(L/K)=1$.  This
property is preserved under finite algebraic extensions, see e.g. (18.1)
in \cite{End}. Note that a Henselian valued field which is defectless
has no proper finite immediate extension. We will use the following
results:

\begin{fact}\vlabel{defectless}Let $(K,v)$ be a valued field. Then $(K,v)$ is defectless
  in the following cases:
\begin{enumerate}
\item $k_K$ has characteristic $0$ (Obvious since then $p=1$).
\item $\Gamma(K)\simeq \zee$ (Theorem 3.3.5 in \cite{EP}).
\item $K$ is finitely generated over a subfield $E$ on which $E$ is
  trivial, and $\trdeg(K/E)=\trdeg (k_K/E)+{\rm
    dim}_\rat\Gamma(K)\otimes \rat$ (Theorem 1.1. in \cite{Ku}). 
\end{enumerate}
\end{fact}

\begin{lem}\vlabel{ram-2} Let $(K,v)$ be a valued field,  $L$ a finite
  Galois extension of 
  $K$, and $M$ an  algebraic extension of $K$. We assume that $K$
  is defectless. 
\begin{enumerate}
\item Assume that $char(K)=p>0$, and that $L$ is an Artin-Schreier
  extension of $K$, generated by a root $\alpha$ of $X^p-X-a=0$. If
  $r(L/K)=p$, then for any $d\in K$, $v(d^p-d-a)<0$. 
\item
Assume that  $M$ is  purely inseparable and finite over $K$.  Then
$r(L/K)=r(LM/M)$, where $M$ is endowed with the unique extension of $v$
to $M$.  
\item Assume that $M$ is Galois over $K$, and that $e(M/K)=1$.  Then
  $e(L/K)$ divides $e(LM/M)$, and equality holds  if the residue fields
  $k_L$ and 
  $k_M$ are linearly disjoint over $k_K$. 
 
\item Assume that $M$ is finite Galois over $K$, and that $e(M/K)$ and $e(L/K)$
  are relatively prime. Then
  $e(L/K)$ divides $e(LM/M)$, and equality holds  if the residue fields
  $k_L$ and 
  $k_M$ are linearly disjoint over $k_K$. 

\item Assume that $M$ is Galois over $K$, that $k_M\subseteq k_K^s$, and that $m \Gamma(L), m\Gamma(M)\subseteq
  \Gamma(K)$ for some integer $m$ relatively prime to $p$. Then also $ m\Gamma({LM})\subseteq \Gamma(K)$. In
  particular, if $m \Gamma({M})=\Gamma(K)$, then $e(LM/M)=1$. 
\item Assume that $M$ is Galois over $K$, that $k_M\subseteq k_K^s$, and
  that $m \Gamma({M})=\Gamma(K)$ for some integer $m$ prime to $p$ and
  divisible by the exponent of the 
  prime-to-$p$ part of $\Gamma(L)/\Gamma(K)$. Then
  $e(LM/M)$ is a power of $p$.   
\end{enumerate}
\end{lem}

\prf % (1) It suffices to show the result for $M=K^{1/p}$. Let $w$ be an extension of
% $v$ to $L$ which ramifies. Note that
% $ML=L^{1/p}$ because $L$ is separable over $K$.  Hence the Frobenius map
% $ML\to L$ sends $M$ to $K$, and if we define 
% $w'$ on $ML$  by $w'(x)=w(x^p)/p$, then it sends $w'$ to $w$. The result
% follows.  \\[0.05in]
(1) Assume that there is $d\in
K$ such that $v(d^p-d+a)>0$. As 
$\alpha+d$ is a root of $X^p-X-(d^p-d+a)$,  we may  assume
that $v(a)>0$. Then $\bar a=0$,  the equation
$X^p-X=0$ has $p$ distinct roots in the residue field, and
therefore $\alpha$ lies in the henselization of $K$. This implies
$e(L/K)=f(L/K)=1$, which contradicts $r(L/K)=p$. Hence, for every $d\in
K$, $v(d^p-d-a)\leq 0$. \\
Assume now that there is  $d\in K$ such that 
$v(d^p-d+a)=0$. As above, we may assume $v(a)=0$. Our  assumption  then
implies that the polynomial $X^p-X-\bar a$ has no root in $k_K$ (since
if $d\in\calo_K$ lifts a root, then $v(d^p-d-a)>0$) and
therefore $\bar\alpha$ generates an extension of degree $p$ of
$k_K$. I.e., $k_L$ is separable
of degree $p$ over $k_K$,  $r(L/K)=1$, and again we reach a
contradiction. This gives us the desired result. \\[0.1in]
For (2) -- (6), we  replace $K$ by its Henselisation, and
therefore assume that $v$ has a unique extension to $K^{alg}$, which we
will also denote $v$. This does
not affect the ramification indices, nor the possible linear
disjointness of the residue fields, nor the defectlessness.    \\[0.05in]
(2) References are to \cite{EP}, chapter 5, and particularly sections
5.2 and 5.3. Reasoning by 
induction, we may assume
that $[M:K]=p$. If $L_v$ denotes the subfield of $L$ fixed by $H=\{\rho\in
\gal(L/K)\mid \forall a\in L\, v(\rho(a)-a)>v(a)\}$, then the residue
field of $L_v$ equals $k_L\cap k_K^s$, 
 $e(L_v/K)$ is
prime to $p$, and $H$ is a $p$-group, $r(L/K)=e(L_v/K)[L:L_v]$. Since $e$ is multiplicative in
towers, and $[M:K]=p$, one clearly has that
$e(L_vM/M)=e(L_v/K)$, and furthermore, the residue field of $L_v$ is a
separable extension of $k_K$, so does not contribute to $r(L/K)$. 
We may therefore assume that $K=L_v$, and the proof is by induction on
$[L:K]$. As $\gal(L/K)$ is a $p$-group, $L$ contains an Artin-Schreier
extension, generated over $K$ by a root $\alpha$ of $X^p-X=a$, for
some $a\in K$. Because $K$ and $M$ are  Henselian and defectless, we know that
$K(\alpha)$ and $M(\alpha)$ are not immediate. Using the multiplicativity in towers of $r$, we will 
therefore assume that  $L=K(\alpha)$. Then $[LM:M]=p$, so that either
$e(LM/M)=p$ or 
$f(LM/M)=p$. It therefore suffices to show that if $e(LM/M)=1$, then the
residue field  $k_{LM}$ of $LM$ is purely inseparable over
$k_K$. \\[0.05in] 
By (1), we know that for all $c\in K$, $v(c^p-c+a)<0$. 
Assume by way of contradiction that $e(LM/M)=1$, and that $k_{LM}/k_K$
is separable. Then $k_{LM}$ is an Artin-Schreier extension of $k_M$,
i.e., is generated over $k_M$ by a root of a polynomial $X^p-X-\bar b$,
for some $b\in M$. Since $LM$ is Henselian, $LM$ contains therefore a
root of $X^p-X-b$, and this root generates $LM$ over $M$. By the theory
of Artin-Schreier extensions\footnote{See e.g. Theorem VIII.8.15 in
  S. Lang, {\em Algebra}, Addison-Wesley 1971.}, there is some positive integer $i<p$ and
element $d\in M$ such that $a+d^p-d=bi$, and in particular,
$v(d^p-d+a)=0$.  Since
$[M:K]=p$, $d^p=c\in K$. But as $L$ is separable over $K$, we have
$K(\alpha)=K(\alpha^p)$ and $\alpha^p$ is a root of
$X^p-X-a^p=0$. Hence,  $v(c^p-c+a^p)=0$, and $L$ is generated over $K$
by a root of 
$X^p-X-(c^p-c+a^p)$. This contradicts (1).  \\[0.05in]
(4) The residue field of $LM$ contains the residue fields of $L$ and
of $M$. Hence $f(LM/K)$ divides $f(L/K)f(M/K)$ and equality holds if
$k_L$ and $k_M$ are linearly disjoint over $k_K$. Each of $e(L/K)$,
$e(M/K)$ divides $e(LM/K)$, and therefore so does their product since
they are relatively prime.  This
proves the first assertion, and the second follows from
$f(LM/K)=f(L/K)f(M/K)$ and the defectlessness and Henselianity of
$K$ (which gives $[LM:K]=e(LM/K)f(LM/K)$). This shows (4),  
 from which also (3) follows. \\[0.05in]
(5) This is implicit in the description of the second exact sequence in
\cite{EP} (section 5.3, see in particular p.~129, and 5.3.3, 5.3.8), but we will give the proof.   Recall that we assume $K$ Henselian, so that $v$ has a unique
extension to $K^{alg}$, also denoted by $v$, and that $p$ denotes the
characteristic of $k_K$ if it is positive, and $1$ otherwise. We let $G=\gal(K^s/K)$,
$G^t=\{\rho\in G\mid 
v(\rho(x)-x)>0\, 
\forall x\in\calo_{K^s}\}$ (the inertia group), $G^v=\{\rho\in G\mid v(\rho(x)-x)>v(x) \,
\forall x\in (K^s)^\times\}$ (the ramification group), and $K^t$, $K^v$, the subfields of $K^s$
fixed by $G^t$ and $G^v$. Then $G^t$ and $G^v$ are normal subgroups of
$G$, $G_v$ is a pro-$p$-group, and $G^t/G^v$ is abelian, all of its
finite quotients having order relatively prime to $p$. The profinite group
$G^t/G^v$ has the following description: Let $\Omega$ denote the
subgroup of $\rat/\zee$ consisting of elements of order prime to $p$. Then
$$G^t/G^v\simeq \Hom(\Gamma(K^s)/\Gamma(K), \Omega).$$
% In particular, any finite quotient of $G^t/G^v$ has order prime to
% $p$. \\[0.05in] 
Note that our assumption on $M$ implies that $M\subseteq K^v$, so that $LM\cap
K^v=(L\cap K^v)M$, and $[LM:LM\cap K^v]=[L:L\cap K^v]$. If $L\not\subseteq K^v$, then $p>1$, and our assumption
that $e(L/K)$ is prime to $p$ implies that  
$$e(L/L\cap K^v)=1, \ \
[L:L\cap K^v]=[LK^v:K^v]=r(L/K)/e(L/K)=[k_L:k_K]_i=[k_{K^vL}:k_{K^v}].$$
Hence $[LM:LM\cap K^v]=[L:L\cap K^v]=[k_L:k_K]_i$, and  $e(LM/LM\cap
K^v)=1$. From this we deduce that $\Gamma(LM)=\Gamma(LM\cap
K^v)$.\\[0.05in] 
Let $G(m)$ denote the subgroup of $G^t/G^v$ corresponding to  those
$f\in \Hom(\Gamma(K^s)/\Gamma(K), \Omega)$ such that $\Ker(f)$ contains
all elements of order $m$ of $\Gamma(K^s)/\Gamma(K)$, and let $K(m)$ be
the subfield of $K^v$ fixed by $G(m)$. Note that $G^t/G(m)$ has exponent
$m$, and is the largest quotient of $G^t$ with this property. In
other words, if $N$ is a Galois extension of $K$, then
$\Gamma(N)/\Gamma(K)$ has exponent dividing $m$ if and only if $N\subseteq
K(m)$. Our
hypothesis on $L$ and $M$ then 
implies that $M\subseteq K(m)$ and $L\cap K^v\subseteq K(m)$. Hence, $LM\cap
K^v\subseteq K(m)$, which implies that $m\Gamma(LM)\subseteq
\Gamma(K)$. This shows the first assertion, and the second follows
immediately:  $\Gamma(LM)=\Gamma(M)=1/m\Gamma(K)$, whence
$e(LM/M)=1$. \\[0.1in]
(6) Let $L_0\subseteq L$ contain $K$ and be such that $[L:L_0]=r(L/L_0)$
equals the highest power of $p$ dividing $r(L/K)$. Then $L_0$ is Galois
over $K$,
 and  $\Gamma(L_0)/\Gamma(K)$ has exponent dividing $m$, i.e., $m\Gamma(L_0)\subseteq \Gamma(K)$. By (5), we get
$e(L_0M/M)=1$. As $[LM:L_0M]$ is a power of $p$, so is  $e(LM/M)=e(LM/L_0M)$.

\para\vlabel{val1}{\bf Generalised power series}. Recall that if $\Gamma$ is an
ordered abelian group and $E$ a field, then the field of generalised
power series $E((t^\Gamma))$ is defined as the set of all
formal sums $f=\sum_{\gamma\in \Gamma}a_\gamma t^\gamma$ with
$a_\gamma\in E$ and such that the support of $f$, 
$\Supp(f)=\{\gamma\in\Gamma\mid a_\gamma\neq 0\}$,  is
well-ordered.\\[0.05in] 
We define a valuation $v$ on $E((t^\Gamma))$ by $v(f)=\inf \Supp(f)$.  If $a\in E((t^\Gamma))$ and $\gamma\in \Gamma$,
we denote by $a\rest\gamma$ the unique element $b$ of $E((t^\Gamma))$
with support contained in $(-\infty, \gamma)$ and such that
$v(a-b)\geq \gamma$. \\[0.05in] 
We denote by $E[t^\Gamma]$, [resp. $E(t^\Gamma)$]  the subring
[resp. subfield] of $E((t^\Gamma))$ generated
by $E$ and all elements $t^\gamma$, $\gamma\in \Gamma$.

\begin{lem}\vlabel{AS} Let $E$ be an algebraically closed field of
  positive characteristic $p$, and let $\Gamma$ be a torsion free
  abelian group. Let $f=\sum a_\gamma t^\gamma\in E[t^\Gamma]$, and let
  $\Delta$ be the
  subgroup of $\Gamma$ generated by ${\rm Supp}(f)$. 
\begin{enumerate}
\item There is $h\in E(t^\Gamma)$ such that $h^p-h=f$ if and only if there
  is a partition of ${\rm Supp}(f)\setminus \{0\}$ into subsets $S_i$
  satisfying the following condition for each $i$: there is $\gamma_i\in
  S_i$ such that for all $\gamma\in S_i$,  
 there exists
  $m(\gamma)\in\nat$ such that 
  $$\gamma=p^{m(\gamma)}\gamma_i,\hbox{ and } \sum_{\gamma\in
    S_i}a_\gamma^{1/p^{m(\gamma)}}=0.$$
Furthermore, if $f=h^p-h$ then $h\in E[t^\Delta]$. 
\item Assume that no element of ${\rm Supp}(f)$ is divisible by $p$ in
  the subgroup 
  $\Delta$, and that for
  some $g\in E[t^\Delta]$, with $|{\rm Supp}(g)|=|{\rm Supp}(f)|$, there is
  some $h\in E[t^\Gamma]$ such that $h^p-h=f-g$. Then there is a (unique)
  bijection $\pi: {\rm Supp}(f)\to {\rm Supp}(g)$ such that for each
  $\gamma\in {\rm Supp}(f)$, there is some $m(\gamma)$ such that
  $\pi(\gamma)=p^{m(\gamma)}\gamma$, and the coefficient  of
  $t^{\pi(\gamma)}$ in $g$ equals $-a_\gamma^{p^{m(\gamma)}}$. 
 
\end{enumerate}
\end{lem}

\prf (1) The right to left implication is clear, since for any integer
$m>0$, 
$x^{p^m}-x=h^p-h$, where $h=\sum_{i=0}^{m-1}x^{p^i}$. Note that this $h$
belongs to the ring generated by $x$. \\
Assume  now that $f=h^p-h$. The ring $E[t^\Gamma]$ is integrally closed,
and therefore $h\in E[t^\Gamma]$. Write $h=\sum_\gamma c_\gamma t^\gamma\in
E[t^\Gamma]$, and assume that the right hand side of the equivalence does not
hold. We will show that this leads to a contradiction. Observe that the
sets $S_i$ are uniquely determined by the support of $f$. 
Using the remark made in the proof of the right to left implication and the fact that $E$ is
perfect, we may assume that no non-zero element of the support of $f$ is
divisible by $p$ in $\Delta$. Hence, each $S_i$ consists of a singleton, and at least one
of them is non-empty. But on the other
hand, we have 
$$a_\gamma=(c_{\gamma/p})^p-c_\gamma\eqno{(*)}$$
for each $\gamma\in\Gamma$. If $\gamma=0$, this simply says that
$a_0=c^p_0-c_0$. If $\gamma\neq 0$,  our condition on the $S_i$'s
implies that $a_{p^i\gamma}=0$ for all
$0\neq i\in\zee$. Hence,  $$c_{p^{i-1}\gamma}^p=c_{p^i\gamma}\hbox{ for
  all }0\neq i\in\zee.$$ But $h$ has finite
support, so this implies that all $c_{p^i\gamma}, i\in\zee,$ are $0$, and
therefore that $a_\gamma=0$, the desired contradiction. \\
Using $(*)$, the same reasoning shows that if $h^p-h=f$, then $h\in E[t^\Delta]$. \\[0.05in] 
(2) This is essentially clear from (1). Note that $0\notin {\rm Supp}(f)$. Let
$S={\rm Supp}(f)\cup {\rm Supp}(g)$, and $\{S_i\}$ the partition of $S$
given by (1). Then each $S_i$ contains at most one element of ${\rm
  Supp}(f)$, has at least two elements (by (1)), and therefore contains
precisely one element of ${\rm Supp}(f)$ and one element of ${\rm
  Supp}(g)$. This gives the bijection $\pi$, and  the rest
follows from (1).

\bigskip\noindent
The following result is probably well-known, but for lack of a reference we
will give its proof.

\begin{prop}\vlabel{val2} Let $E$ be a field, and
  $\Gamma$  a 
subgroup of $\ree$ such that $\bigcap_{n\in\nat}p^n\Gamma=(0)$ if
$char(E)=p>0$. 
The elements of $E((t^\Gamma))$ which are
separably algebraic over $E(t^\Gamma)$ have the following property: their
support is either finite or of order type $\omega$ and cofinal in
$\ree$. Moreover, the subring $E[t^\Gamma]$ is
dense in $E(t^\Gamma)^s\cap E((t^\Gamma))$. \end{prop}

\prf We consider the completion $E(t^\Gamma)^c$ of $E(t^\Gamma)$ with
respect to the valuation. This is  the smallest field
containing limits of all sequences $(a_n)_{n\in\nat}$ of elements of
$E(t^\Gamma)$, with
$v(a_{n+1}-a_n)$  increasing and cofinal in $\ree$. Then
$E(t^\Gamma)^c$ is henselian (see Chapter 2 of \cite{[S]}).

\smallskip\noindent
{\bf Claim}. $E(t^\Gamma)^c$ coincides with the set $K$ of elements of
$E((t^\Gamma))$ whose support is either finite or of order type
$\omega$ and cofinal in $\Gamma$.\\[0.05in]
Note that $a\in K$ if and only if for every
$\gamma\in\ree$, $\Supp(a)\cap (-\infty,\gamma)$ is finite. Clearly $K$
contains $E[t^\Gamma]$ and is closed under addition. For
multiplication, let $a,b\in K$, $\gamma\in \ree$; then
$\Supp(ab)\cap (-\infty,\gamma)\subset \Supp(a\rest{\gamma-v(b)}\cdot
b\rest{\gamma-v(a)})$. Similarly, if $a\in K$ with $v(a)>0$,
$\gamma>0$,
and $k$ is the smallest integer such that $kv(a)>\gamma$, then
$\Supp((1+a)\inv)\cap (-\infty,\gamma)\subset
\Supp(1+a\rest\gamma+\cdots+a^{k-1}\rest\gamma)$. This implies that $K$
is a subfield of $E((t^\Gamma))$. It follows that every element of
$E(t^\Gamma)^c$ is the limit of a sequence of elements of
$E[t^\Gamma]$, and therefore coincides with $K$. This finishes the proof
of the claim. 

\medskip\noindent
Thus the ring
$E[t^\Gamma]$ is dense in $E(t^\Gamma)^c$:
for every $a\in E(t^\Gamma)^c$ and $\gamma\in \ree$, there is $b\in
E[t^\Gamma]$ such that $v(a-b)>\gamma$. 
The field $E(t^\Gamma)^c$ is henselian, and if $char(E)=0$, it 
has no proper
immediate algebraic extension (see \ref{defectless}), so that the result holds. 

\medskip\noindent
Assume now that
$char(E)=p>0$, and that $L$ is a finite Galois extension of
$E(t^\Gamma)^c$, with Galois group $G$.
We need to show that $L_0:=L\cap E((t^\Gamma))= E(t^\Gamma)^c$. We
assume this is not the case, and take such an $L_0$ of minimal degree
over $K$, and $L$ its normal closure over $K$. We
denote also by $v$ the unique extension of $v$ to $L$. If $G_v=\{\rho\in
G\mid \forall x\in L\, v(\rho(x)-x)>v(x)\}$, then the subfield $L_1$ of $L$
fixed by $G_v$ is contained in  $L_2:=E'(t^{\Gamma'})^c$, for some separable
extension $E'$ of $E$ and  overgroup
$\Gamma'$ of $\Gamma$  such that
$[\Gamma':\Gamma]$ is prime to $p$. 
Indeed, $L_1=K(\alpha,\beta)$, where the residue field of $K(\alpha)$ is
separable over $E$, and the extension $K(\alpha,\beta)/K(\alpha)$ is
totally ramified of degree prime to $p$, see Corollary 5.3.8 in
\cite{EP}. Moreover 
$L_2$ and $L_0$ are linearly disjoint over $K$, and $\gal(L/L_1)=G_v$ is a
$p$-group, which is normal in $G$. The minimality of $L_0$ then implies that
$L_0L_1$ contains a Galois extension of degree $p$ of $L_1$, i.e.,
contains some 
Artin-Schreier extension of $L_1$ which is immediate over $L_1$. As
$L_2$ and $L_0$ are linearly disjoint over $K$, it is enough to show: no Artin-Schreier extension of $L_2$
is immediate. \\[0.05in] 
Assume by way of contradiction that $L_2(a)$ is an immediate
Artin-Schreier extension of $L_2$, where $b=a^p-a\in L_2$. 
 Because $L_2$ is Henselian, every
element of positive valuation is of the form $x^p-x$ for some $x\in L_2$, so that  $v(b)\leq
0$, and without loss of generality, $b\in E'[t^{\Gamma'\leq 0}]$ 
has finite support (use the claim, and the fact that every element of
positive valuation is of the form $x^p-x$ in a Henselian field). Write $b=\sum_\gamma c_\gamma t^\gamma$,
$c_\gamma\in E'$, $\gamma\in \Gamma'$.
%\sum_{i=1}^n c_i t^{\gamma_i}$, with
%$\gamma_1<\gamma_2<\cdots <\gamma_n\leq 0$. 
We may then assume that each
$c_\gamma t^\gamma$ with $\gamma<0$ is not a $p$-th power in $L_2$:
otherwise, let $s$ be 
maximal such that $c_\gamma t^\gamma\in L_2^{p^s}$ (recall that no element
of $\Gamma$ is infinitely $p$-divisible, hence the same holds for
$\Gamma'$), and add to $b$ the 
element $c_\gamma^{p^{-s}}t^{p^{-s}\gamma}-c_\gamma t^\gamma$ (which is
of the form $d^p-d$ for some $d\in L_2$, see \ref{AS}, proof of (1)). Let
$\gamma_0$ be the least element of the support of $b$. Now,
$v(a)=v(b)/p$, and because $e(L_2(a)/L_2)=1$, it follows that $\gamma_0\in
p\Gamma'$; if $\gamma_0<0$, it then implies that $c_{\gamma_0}$ is not a
$p$-th power in $E'$, but  as $v(a
-c_{\gamma_0}t^{\gamma_0/p})>v(a)$, we get a contradiction with
$f(L_2(a)/L_2)=1$. Finally, we cannot have $\gamma_0=0$, since then
$L_2(a)$ would be a
purely residual extension of $L_2$. \\[0.01in]
{\bf Comment}. Actually, this proof can be easily modified to show that
$K=E(t^\Gamma)^c$ is defectless, the point being to reduce to the case of immediate
Artin-Schreier extensions.

\para\vlabel{alg1}{\bf Conventions and some basic results and definitions in
algebraic geometry}. Recall that we work in a large saturated model
$(\Omega,\si)$ of ACFA. Our varieties will always be quasi-projective and
absolutely irreducible. If a variety $V$ is defined over the field $K$,
then $K(V)$ denotes its function field, and if $V$ is affine, $K[V]$ 
denotes the affine ring of $V$ (over $K$). We view the elements of
$K(V)$ as (partially defined) functions on $V$. % The 
% function field of the variety $V$ will be denoted b
%y $\Omega(V)$,%
%and viewed as a set of rational functions on $V$.
If $x\in V(K)$, then $\calo_{x,V}$ will denote the ring of elements of
$K(V)$ which are defined at $x$, and $\calm_{x,V}$ its maximal ideal,
which consists of functions vanishing at $x$. We use the same notation
if $U$ is an affine open subset of $V$, and $x$ a prime ideal of $K[U]$,
to denote the localization of $K[U]$ at the ideal $x$ and its maximal
ideal. \\[0.1in]
%
%
%{\bf Normalisation}. 
Recall that the variety $V$ is {\em normal} if whenever $U$ is an
affine open subset of $V$ (defined over $K$), then $K[U]$ is
integrally closed in its fraction field. If $V$ is normal and affine,
and $L$ is a 
finite algebraic extension of $K(V)$, then the {\em normalisation of $V$
  in $L$} is  the variety  $W$ whose affine ring is the integral closure $R$
of $K[V]$ in $L$; it is therefore affine, and comes with a dominant finite map $f:W\to V$, dual
to the inclusion $K[V]\subset R$. 
If $V$ is normal and projective, and $V=\bigcup U_i$
where the $U_i$'s are affine open, then the {\em normalisation $W$ of $V$ in
$L$} is $W=\bigcup W_i$, where each $W_i$ is the normalisation of $U_i$
in $L$. The map $f:W\to V$ is obtained by glueing the maps $f_i:W_i\to
U_i$.  See Theorem 4 of III.8 in \cite{Mu}, or Theorem 4 of V.4 in
\cite{La72}.

% A dominant morphism $f:V\to W$ induces a dual morphism
% $f^*:\Omega(W)\to \Omega(V)$.

\para \vlabel{ram-1}{\bf Finite morphisms, ramification.} A dominant
morphism  $f:V\to W$ between two affine varieties is {\em finite} if
$K[V]$ is integral over $f^*(K[W])$. If $V$ and $W$ are 
quasi-projective, then it is {\em finite} if $W$ can be covered by
affine subsets $W_i$ with $f\inv (W_i)$ affine, and the restriction of
$f$ to $f\inv (W_i)$ finite. \\
Assume that $W$ is {\bf normal}, and $f$ finite. If $y\in W(\Omega)$, the dominant finite morphism $f$ is {\em unramified 
over $y$} if  $f\inv(y)$ contains
$\deg(f)$ distinct points, where $\deg(f)=[K(V):f^*K(W)]$; if
$|f\inv(y)|< \deg(f)$, then we say that $f$ {\em ramifies over} $y$. The set of
points of $W$ over which $f$ is ramified is % Zariski closed and is
called the {\em locus of
ramification} (or {\em ramification locus}) of $f$. %  is the Zariski
% closure of the set of ramification points of $f$.
 This is a Zariski closed subset of $W$, which is defined over $K$ (see Thm 4 of II. 6.3 of \cite{Sh1}). If
$f$ is inseparable, it is all of $W$. We say that $f$ is {\em unramified
over $W$}, if it is unramified over every point of $W$. 

There is an alternate,
equivalent definition involving local rings, see \cite{Mi} I.3 (p. 21): without loss of generality
assume that $V$ and $W$ are affine, let $x$ be a prime ideal of $K[V]$,
$y=f(x)$, and $\calo_x,\calo_y$ the corresponding localisations of
$K[V]$ and $K[W]$; we identify $\calo_y$
with a subring of $\calo_x$ via $f^*$. Then $f$ is
unramified at $x$ if and only if $\calo_x/\calm_y\calo_x$ is a separable
field extension of $\calo_y/\calm_y$. Note that this statement has two
implications: that $\calm_y$ generates the maximal ideal of $\calo_x$, and
that  the residue field $\calo_x/\calm_x$ is separably algebraic over
$\calo_y/\calm_y$. We say that $f$ is unramified over $y$ if it is
unramified at every point $x\in f\inv(y)$, and that it is unramified
over $W$ if it is unramified over every point of $W$ (or equivalently, of
$W(K^{alg})$).  %\\[0.1in]
 % If $f:V\to W$ is finite unramified
% and flat, then $f$ is {\it \'etale}.
%
%
\para\vlabel{ram-1bis} {\bf Ramification divisor}. Let $f:V\to W$ be as above, and $S\subset W$ a subvariety of
codimension $1$. As $W$ is normal, the function which to an element
$g\in K(W)$ associates its order of vanishing on $S$, defines a discrete
valuation $v_S$ on $K(W)$, which we call a {\em divisorial
  valuation}. This valuation $v_S$ extends to $K(V)$, and may or may not
ramify. Note that by \ref{defectless}, $v_S$ is defectless. \\
The {\em ramification divisor} is the union of all subvarieties $S$ of $V$
of codimension $1$ such that the associated valuation $v_S$ ramifies in
$K(W)$. {Note that we are not interested in multiplicities here; thus the ramification divisor in our sense
is the support of the ramification divisor in the sense, e.g., of
Hartshorne \cite{Ha}, p. 301 (for curves).}  \\[0.1in]
We will use the following result of Zariski, see \cite{L56}, Proposition 4, or
\cite{Ab}, Thm 1:
\begin{thm}\vlabel{Zar}
{ If $W$ is non-singular, then the ramification divisor of $f$ and the
ramification locus of $f$ coincide.} \end{thm}

\begin{lem}\vlabel{su1} Let $E=\acl(E)$,
let $B$ be a definable  subgroup
of some algebraic group $G$ defined over $E$, and let
$a$ be a generic of $B$.  Assume that  $\si(a)\in E(a)^{alg}$, let
$\ell>0$ and let $B_{(0)}$ and $V_{0,\ell}$ be the algebraic loci of $a$ and
of 
$(a,\si^\ell(a))$ 
over $E$.
Assume that  the variety $V_{0,\ell}$ contains a proper infinite
subvariety
$W$ such that
$\pi_0(W)^{\si^\ell}=\pi_\ell(W)$ have the same dimension as $W$, where
$\pi_0:V_{0,\ell}\to B_{(0)}$ and $\pi_\ell:V_{0,\ell}\to
B_{(0)}^{\si^\ell}$ are the natural projections given by the inclusion $V_{0,\ell}\subset
B_{(0)}\times
 B_{(0)}^{\si^\ell}$.  Then $\evSU(a/E)>1$.
\end{lem}

\prf The assumption $\si(a)\in
E(a)^{alg}$ implies that $\dim (V_{0,\ell})=\dim(B_{(0)})$, and that
$\dim(W)<\dim(B_{(0)})$.  Moreover,
$\si(a)\in E(a)^{alg}$ implies that $a$ is also a generic of $B(\ell)$,
and therefore $B$, $B(\ell)$ and $tp(a/E)$ have the same
evSU-rank. \\[0.05in]  
We now work in $\Omega[\ell]$. Let $F=\acl(F)$ contain $E$ and such that
$W$ is defined over $F$. By properties of the theory ACFA and the
fact that $\Omega[\ell]$ is a model of ACFA, there is some $b\in
B(\ell)$ such that $(b,\si^\ell(b))\in W$, and
$\trdeg(F(b)_\si/F)=\dim (W)$. Then
$0<\trdeg(F(b)_\si<\trdeg(E(a)_\si/E)$, so $b$ is not a generic of
$B(\ell)$, and is not algebraic over $F$. Hence $\SU(a/E)[\ell]>1$. 
% such
% that $b$ is a generic of $$  and $(b,\si^\ell(b))\in W$. 
% By Theorem \ref{f1}, the set
% defined by $\{x\in G(\Omega)\mid (x,\si^\ell(a))\in W\}$ is a Boolean
% combination of cosets of quantifier-free definable subgroups of the
% connected component $B(\ell)^0$ of $B(\ell)$, and these subgroups must
% be proper as $\dim(W)<dim(B_{(0)})$. I.e., $B(\ell)$
% contains some 
% definable subgroup $C$ which is infinite and of infinite index in
% $B(\ell)$.
% Thus $\SU(a/E)[\ell]>1$.

\begin{prop}\vlabel{ram0} Let $E=\acl(E)$, $G$ a connected commutative
  algebraic group, and 
  $B$ a definable modular subgroup of  $G(\Omega)$, of eventual SU-rank
  $1$ and Zariski dense in $G$,
  everything being defined over $E$. Let $a$ be a
  generic of $B$ over $E$ (so that $E(a)\simeq E(G)$), and assume that $\si(a)\in E(a)^{alg}$. % If
%   $\trdeg(E(a)/E)=1$, we assume in addition that $E(a)_\si\subseteq
%   E(a)^s$. 
  Let $L$ be a finite Galois extension of $E(a)$, and assume
  that $L$ is 
linearly
  disjoint from $E(a)_\si$ over $E(a)$, and is such that
  $L(\si(a))=\si(L)(a)$.  If  $W$ is 
  the normalisation of $G$ in $L$, then the map $f:W\to G$ is
  unramified. 
\end{prop}

\prf Assume by way of contradiction that the map $f: W\to G$ is
ramified, let $\cals$ be the ramification divisor of $f$,  
 $S$  an irreducible component of $\cals$, and $v=v_S$ the
associated valuation (see
\ref{ram-1bis}). Let $V_1$ be the algebraic locus of $(a,\si(a))$ over
$E$. Then $V_1$ is an algebraic subgroup of $G\times G^\si$, and is
therefore non-singular.  

Each of the projections $\pi_0:V_1\to G$ and $\pi_1:V_1\to
G^\si$, being a 
group homomorphism, is the composition of a purely 
inseparable morphism with an unramified  morphism. Moreover, as $E$ is
algebraically closed, so that $G(E)$ contains the torsion subgroup of $G$,
the field $E(a,\si(a))$ is a normal extension of $E(a)$ and of
$E(\si(a))$. Then, Lemma \ref{ram-2} (2) and (3) imply that if $w$ is
any extension of $v$ to $E(a,\si(a))$, then
$r(L/E(a))=r(L(\si(a))/E(a,\si(a)))=r(\si(L)/E(\si(a)))$ (here we
use the fact that $\si(L)(a)=L(\si(a))$). Then 
the restriction of $w$ to $E(\si(a))$ is a divisorial valuation $v_T$, 
with $T$ an irreducible component of $\pi_1\pi_0\inv(S)$ (this uses the
fact that the ramification divisor coincides with the ramification
locus -- see \ref{Zar} -- and that the maps $\pi_0$ and $\pi_1$ are
finite, everywhere defined).  

On the other
hand, we know that $\cals^\si$ is the ramification divisor of  $f^\si:
W^\si\to G^\si$. Hence $T$ is an irreducible component of
$\cals^\si$. This reasoning can be repeated: we fix an extension $w$ of
$v$ to $E(a,\si(a),\ldots,\si^\ell(a))$, with $\ell= |\cals|$. We
then get a sequence $S=S_0, S_1,\ldots,S_\ell$ of elements of $\cals$ such
that the restriction of $w$ to $E(\si^i(a))$ coincides with the
valuation $v_i$ associated to 
${S_i^{\si^i}}$. Then necessarily there are $0\leq i<j\leq \ell$ such that
$S_i=S_j$, and looking at the algebraic locus $V_{i,j}$ of $(\si^i(a),\si^j(a))$
over $E$, we obtain a subvariety $U$ of $V_{i,j}$, such that $v_U=w\rest{E(\si^i(a),\si^j(a))}$
ramifies in $\si^i(L)(\si^j(a))$, and the projections of $U$ on
$G^{\si^i}$ equals $S_i^{\si^i}$, the projection of $U$ on $G^{\si^j}$
equals $S_i^{\si^j}$. Lemma \ref{su1} gives us the desired contradiction
when $\dim(G)>1$. 

Assume now that $\dim(G)=1$.  Because $tp(a/E)$ is modular and ${\rm tr.deg}(a/E)=1$,
 we know that the set $\{[E(a,\si^k(a)):E(a)]_s\mid k\in\zee\}$ is
 unbounded (see (4.5) in \cite{[CH]}). Choose $k\in\zee$ such that
 $[E(a,\si^k(a)):E(a)]_s=N>|\cals|$.   As  $E$ is
algebraically closed, if  $P\in \cals$ then the valuation
 $v_P$ has $N$ distinct extensions $w_1,\ldots,w_N$ to $E(a,\si^k(a))$,
 which correspond 
 to the $N$ distinct points $Q_1,\ldots,Q_N\in G^{\si^k}(E)$ such that
each $(P,Q_i)$ belongs to the algebraic locus $V_{0,k}$ of
$(a,\si^k(a))$ over $E$; the restrictions  of the valuations $w_i$ to 
 $E(\si^k(a))$ are therefore distinct (they equal the $v_{Q_i}$), and
 the points $Q_1,\ldots,Q_N$ are in $\cals^{\si^k}$. This gives  the
desired contradiction.

\begin{lem}\vlabel{lem1}  Let $G$ be a semi-abelian variety
defined over $E=\acl(E)$, and let $B$ be a definable subgroup of
$G(\Omega)$ of finite SU-rank. Let $n\in\nat^{>0}$, and let $a$ be a generic of
$B$ over $E$, and $b$ such that $[n]b=a$. Then $E(b)_\si$ is a finite
$\si$-stable normal extension of $E(a)_\si$.\end{lem}

\prf We use the notation introduced at  the beginning of section 2. The truth of this statement only depends on $qftp(a/E)$ (see
(2.9)(2) in \cite{[CH]}), and we
may therefore assume that $b\in B$, and that $B=\tilde B$. Since
$[B:B^0]$ is finite and $E$ is algebraically closed, we may also assume
that $B=B^0$. The normality of $E(b)_\si$ over $E(a)_\si$ follows from
the fact that $E$ contains the torsion subgroup of $G$ and of the
$G^{\si^i}$. Let $m$ be such
that $B=\tilde B_{(m)}$, and let $N\geq m$. Then $B_{(N)}$ is a
semi-abelian variety which projects onto $B_{(m)}$ with finite
fibers, so that $\dim(B_{(N)})=\dim(B_{(m)})=: g$. Moreover, if $r$ is the
dimension of the maximal abelian quotient of $B_{(m)}$, then it is also
the dimension of the maximal abelian quotient of $B_{(N)}$. Therefore,
the map $[n]:B_{(N)}\to B_{(N)}$ has degree
$n^{2r+g-r}=n^{g+r}$.\\[0.05in]
The tuple 
$(b,\si(b),\ldots,\si^N(b))$ is a generic
of $B_{(N)}$ for every $N$, whence 
$$[E(b,\si(b),\ldots,\si^N(b)):E(a,\ldots,\si^N(a))]=n^{g+r}, \qquad
\hbox{and }[E(b)_\si:E(a)_\si]=n^{g+r}.$$

\begin{lem}\vlabel{lem31} Let $E=\acl(E)\subseteq F=\acl(F)$,
and let $a$ be independent from $F$ over $E$. Assume that $L$ is a
proper finite separable $\si$-stable extension of $\acl(Ea)F$. Then there is
$F_1$, independent from $(a,F)$ over $E$, and a finite $\si$-stable
Galois 
extension $M$ of $FF_1(a)_\si$ such that $L\subseteq
\acl(F_1a)M$.
\end{lem}

\prf Write $A=E(a)_\si^s$, and $L_0=L\cap (AF)^s$. Note that since
$L$ is separable over $\acl(Ea)F$, and $\acl(Ea)$ is purely inseparable
over $A$, we have 
$L=L_0\acl(Ea)$. It therefore suffices to find such an $M$ with
$L_0\subseteq \acl(F_1a)M$. \\
$AF$ is a normal extension of $F(a)_\si$, and the normal 
closure  of $L_0$ over $F(a)_\si$ is therefore a finite $\si$-stable
Galois extension of $AF$ (see (2.9)(3) in \cite{[CH]}). We may
therefore assume that  $L_0$ is Galois over $F(a)_\si$.  \\
Recall that we work inside the large
difference field $\Omega$. Let $\varphi$ be an $A$-automorphism
of $\Omega$ such that $\varphi(F)=F_1$ is independent from $(F,a)$ over
$E$, and let $L_1=\varphi(L_0)$. Then $L_0L_1$ is a Galois extension of
$FF_1(a)_\si$, which contains $AF$, and we will
identify $$\gal(L_0L_1/FF_1(a)_\si)\simeq
\gal(L_0/F(a)_\si)\times_{\gal(A/E(a)_\si)}\gal(L_1/F_1(a)_\si).$$
Consider the subgroup $H$ of
$\gal(L_0L_1/FF_1(a)_\si)$ defined by
$$H=\{(\tau,\varphi\tau\varphi\inv)\mid \tau\in \gal(L_0/F(a)_\si)\}.$$
As $\varphi$ is an $A$-isomorphism of difference fields,
$H$ is a closed subgroup of $\gal(L_0L_1/FF_1(a)_\si)$, which projects onto
$\gal(A/E(a)_\si)$ and $\si\inv H\si=H$.
Observe that
$[\gal(L_0L_1/FF_1(a)_\si):H]=[L_0:AF]$ is finite. Hence, the
subfield $M$ of $L_0L_1$ which is fixed by $H$ is a finite $\si$-stable
extension of $FF_1(a)_\si$, and intersects $AFF_1$ in
$FF_1(a)_\si$.

We have  $\gal(L_0L_1/L_1)\cap H=\{1\}$, whence
$L_1M=L_0L_1$. From $L_1\subset F_1(a)_\si^s$, we obtain the result.

\begin{lem} \vlabel{lem4} Let $E=\acl(E)$, $a$ a finite tuple,
and assume that $M$ is a finite $\si$-stable  extension of
$E(a)_\si$. Then there are $r\in\nat$ and a finite  extension $L$
of  $E(a,\ldots,\si^r(a))$ such that $M=LE(a)_\si$, $L$ is linearly
disjoint from $E(a)_\si$ over $E(a,\ldots,\si^r(a))$, 
$[L:E(a,\ldots,\si^r(a))]=[M:E(a)_\si]$ and
$\si(L)(a)=L(\si^{r+1}(a))$. \\
Moreover, if $L$ is as above and  $E(a,\ldots,\si^r(a))\subseteq L_0\subseteq L$ is such that $L_0E(a)_\si$ is $\si$-stable,
then we also have $\si(L_0)(a)=L_0(\si^{r+1}(a))$, and $L\si(L)=L_0\si(L_0)=\si(L_0)L$. 
\end{lem}

\prf Fix a finite tuple $\alpha$ generating $M$ over $E(a)_\si$. Then there are
integers $i\leq j$ such that the extension
$M'=E(\si^i(a),\ldots,\si^j(a),\alpha)$ of
$E(\si^i(a),\ldots,\si^j(a))$ has  degree $[M:E(a)_\si]$.
By assumption, $\si(M')\subseteq M'E(a)_\si$ and $M'\subseteq
\si(M')E(a)_\si$. Hence there are integers $\ell\leq i$ and $k\geq j$
such that $\si(M')\subseteq M' E(\si^\ell(a),\ldots,\si^k(a))$ and
$M'\subseteq \si(M')E(\si^\ell(a),\ldots,\si^k(a))$. Define
$L=\si^{-\ell}(M')E(a,\ldots,\si^{-\ell+k}(a))$. Then $[L:E(a,\ldots,
\si^{-\ell+k}(a))]=[M:E(a)_\si]$. 
An easy computation shows that $\si(L)\subseteq L(\si^{-\ell+k+1}(a))$
and $L\subseteq \si(L)(a)$.\\
Write $b=(a,\ldots,\si^{-\ell+k}(a))$, and let now $E(b)\subset
L_0\subset L$ be such that $L_0E(b)_\si$ is $\si$-stable. The linear
disjointness of $L$ and $E(a)_\si=E(b)_\si$ over $E(b)$ implies that $L$
is linearly disjoint from $M_0=L_0E(b)_\si$ over $L_0E(b)=L_0$. In
particular, $[L:L_0]=[L\si(L_0):L_0\si(L_0)]=[L\si(L):L_0\si(L_0)]$
(because $\si(b)\in \si(L_0)$), and therefore we obtain
$[L_0\si(L_0):E(b,\si(b))]=[L_0:E(b)]$, which implies
$L_0(\si(b))=L_0\si(L_0)=\si(L_0)(b)$. The second set of equalities is
also clear from the above. 

\begin{lem}\vlabel{lem5} Let $E=\acl(E)$, $B$ a  $E$-quantifier-free definable
  subgroup of a semi-abelian variety $G$, and assume that if $a\in B$, then
  $\si(a)\in E(a)^{alg}$.  Let $n\geq 1$ be an integer,  let
  $B(n)$ be the  $\si^n$-closure of $B$ (for
  the $\si^n$-topology), let $B^0$ be the connected component of $B$,
  and let $m\in\zee$, $\tau=\frob^m\si^n$. We use the additive notation
  for the group law. The
  following properties are equivalent:
\begin{itemize}
\item[(a)] For all (some) generic $a$ of $B^0$ (over $E$), and field $F=\acl(F)$
  containing $E$, if $L$
  is a finite separable $\si$-stable  extension of $F(a)_\si$, then for
  some integer $N$ and $b\in G(\Omega)$ such that $[N]b=a$, we have $L\subset F(b)_\si$.
\item[(b)] For all $a\in B$ and field $F=\acl(F)$ containing $E$, if $L$
  is a finite separable $\si$-stable  extension of $F(a)_\si$, then for
  some integer $N$ and $b\in G(\Omega)$ such that $[N]b=a$, we have $L\subset F(b)_\si$.
\item[(c)] For all $a\in B(n)$ and field $F={\rm acl}_{\si^n}(F)$ containing $E$, if $L$
  is a finite separable $\si^n$-stable  extension of $F(a)_{\si^n}$, then for
  some integer $N$ and $b\in G(\Omega)$ such that $[N]b=a$, we have $L\subset F(b)_{\si^n}$.
\item[(d)] For all $a\in B(n)$ and field $F={\rm acl}_{\tau}(F)$ containing $E$, if $L$
  is a finite separable  $\tau$-stable  extension of $F(a)_{\tau}$, then for
  some integer $N$ and $b\in G(\Omega)$ such that $[N]b=a$, we have $L\subset F(b)_{\tau}$.
\end{itemize}
\end{lem}

\prf  First some comments about (a): recall that the finite
$\si$-stable Galois extensions of a difference field are
determined by its isomorphism type and do not depend on the extension
of the automorphism to the algebraic closure (\ref{bab}). Hence, the
truth of (a) only depends on the generic quantifier-free
type of $B^0$, 
not on a particular realisation.  
(b) clearly implies (a), and we will now show that (a) implies
(b). 
Let $a\in B$, $F=\acl(F)$ containing $E$ and $L$ a finite separable 
$\si$-stable extension of $F(a)_\si$. Let $a_2$ be a generic of $B^0$
over $\acl(Fa)$, let $a_1=a-a_2$, and $F_1=\acl(Fa_1)$. By (a), there 
are an integer $N$ and a tuple 
$b_2$ such that $[N]b_2=a_2$ and $F_1L\subset F_1(b_2)_\si$. If $b_1\in F_1$ is
such that $[N]b_1=a_1$, setting $b=b_1+b_2$, we obtain $L\subset
F_1(b)_\si\cap \acl(Fa)=F(b)_\si$. This proves (b) and finishes the
proof of the first equivalence. \\[0.05in]
Observe that if $a\in B$ realises a generic of $B$, then 
$qftp(a/E)[m]$ is a generic type of $B(m)$. Together with the first
equivalence (applied to $\si^n$ and to $\tau$), this  allows us to restrict our attention to tuples $a$
which are generics of $B^0$. We let $a$ denote a generic of $B$ over
$E$. \\[0.05in] 
Assume (c), and let $L$ be a finite separable
$\si$-stable extension of $F(a)_\si$, for some $F=\acl(F)$ containing
$E$. Because $\si(a)\in E(a)^{alg}$, it follows that for some $N>0$, we
have $[N]\si^j(a)\in E(a)_{\si^n}$ for all $j=1,\ldots, n-1$, and that $F(a)_\si$ is a finite
$\si^n$-stable extension of $F(a)_{\si^n}$ (Lemma \ref{lem1}). Hence so is $L$, which
easily gives the result: if $L_0=L\cap F(a)_{\si^n}^s$, then $L_0$ is a
finite separable $\si^n$-stable extension of $F(a)_{\si^n}$ and
$L=L_0F(a)_\si$, so that (c) implies (b). \\[0.05in]
Assume now (b). Let $F={\rm
  acl}_{\si^n}(F)$ contain 
$E$, and $L$ a finite separable $\si^n$-stable extension of $F(a)_{\si^n}$. We may assume that $F$ is closed under $\si$ as well
(by (1.12) in \cite{[CH]} applied to the extension $FE(a)_\si$ of
$E(a)_\si$). But then $L\si(L)\cdots
\si^{n-1}(L)$ is a finite separable $\si$-stable extension of $F(a)_\si$, and
therefore contained in $F(b_1)_\si$ for some $b_1$ with $[N]b_1=a$ some
$N$. As $\si(a),\ldots,\si^{n-1}(a)\in E(a)^{alg}$, there is $M$ such
that $E(a)_\si\subset E(b_2)_{\si^n}$ for some $b_2$ with 
$[M]b_2=a$. If $b$ is such that $[MN]b=a$, then $L\subset
E(b)_{\si^n}$. \\ [0.05in] 
To show the last equivalence, we may assume that $n=1$. Observe that
$\tau$ is definable in $(\Omega,\si)$, that $\si$ is definable in
$(\Omega,\tau)$, and that $E(a)_\si E(a)_\tau$ is purely
inseparable over $E(a)_\si$ and over $E(a)_\tau$. \\[0.05in]
Assume first
that it holds for $\tau$, and that $m<0$;  let $F=\acl(F)$ ($={\rm
  acl}_\tau(F)$) contain $E$, and 
 $L$  a finite separable $\si$-stable  extension
of $F(a)_\si$. Apply Lemma \ref{lem4} to get a generator $\alpha$ of
$L$ over $F(a)_\si$, and an integer $r$ such that the minimal
polynomial of $\alpha$ over $F(a)_\si$ has its coefficients in
$F(a,\si(a),\ldots,\si^r(a))$ and
$F(a,\si(a),\ldots,\si^r(a),\si^{r+1}(a),\alpha)=F(a,\si(a),\ldots,\si^{r+1}(a),\si(\alpha))$. Note
that $\tau(\alpha)=\si(\alpha)^{p^{m}}$  generates $\si(L)$ over
$E(\si(a),\ldots,\si^{r+1}(a))$ (because $\si(L)$ is separable over $E(\si(a),\ldots,\si^{r+1}(a))$). Then 
$\si=\tau\frob^{-m}$, so that
$F(a,\si(a),\ldots,\si^r(a),\si^{r+1}(a))\subset F(a,\tau(a),\ldots,
\tau^{r+1}(a))$. As  the extension $F(a,\tau(a),\ldots,
\tau^{r+1}(a))$ of $F(a,\si(a),\ldots,\si^r(a),\si^{r+1}(a))$ is purely
inseparable, it follows that $$F(a,\tau(a),\ldots,
\tau^{r+1}(a),\alpha)=F(a,\tau(a),\ldots,
\tau^{r+1}(a),\tau(\alpha)),$$ so that $F(a)_\tau(\alpha)$ is a finite
 separable $\tau$-stable extension of $F(a)_\tau$. By assumption, there is some $N$
and $b$ such that $[N]b=a$ and $\alpha\in F(b)_\tau$. Our
assumptions on $\alpha$ then imply $\alpha\in F(b)_\si$.\\
If $m=0$, there is nothing to prove, and if $m\geq 1$,  replacing $\tau$
and $\si$ by 
 $\tau\inv$ and $\si\inv$ gives the result.  The other direction  is
 symmetric, interverting the roles of $\si$ and $\tau$. 

\begin{remark} \vlabel{lem6} Let $E=\acl(E)$, $a$ a tuple in $\calu$ and
  assume that $\evSU(a/E)=1$, and $\trdeg(E(a)_{\si^\ell}/E)$ does not
  depend on $\ell>0$. Let $\tau=\frob^m\si^n$, and consider $a$
  in the reduct $(\Omega,\tau)$. Then also $\evSU_\tau(a/E)=1$.
\end{remark}

\prf Our assumption on the transcendence degrees implies that
$E(a)_\si\subset E(a)_{\si^n}^{alg}$. By definition of the $\evSU$-rank,
we know that the evSU-ranks of $a$ over $E$ in $(\Omega,\si)$ and in
$(\Omega,\si^n)$ are the same, and so both equal $1$. Assume that there
is some $F={\rm acl}_\tau(F)$ containing $E$ and such that $a$ is not
independent 
from $F$ over $E$ (in $(\Omega,\tau)$, equivalently in
$\Omega[n]$). Noting that $F={\rm acl}_{\si^n}(F)$, implies $a\in F$ and
gives the desired conclusion.

\para\vlabel{val3}{\bf Algebraically closed valued fields}. Consider the
theory of
algebraically closed valued fields in the $2$-sorted language
$\call_{\rm val}$, with sorts for the valued field and the value group,  $v$  the valuation map $:K\to \Gamma\cup\{\infty\}$,
and $K$ and $\Gamma$ are structures in the languages $\{+,-,\cdot,0,1\}$ and
$\{+,-,0,\leq\}$ respectively. V.~ Weispfenning \cite{[W]} showed that this theory
eliminates quantifiers. (Other quantifier elimination results for
algebraically closed valued fields were obtained by A. Robinson \cite{[R]}
and F.~Delon \cite{[D]}.)\\[0.05in]
Let $(K,v)$ be an algebraically closed valued field, $E$ an algebraically
closed
subfield on which $v$ is trivial, and assume that $\bar a=(a_1,\ldots,
a_n)\in K$ is such that $v(a_1),\ldots,v(a_n)$ are $\rat$-linearly
independent. Then, if $f(X_1,\ldots,X_n)=\sum_{\nu}b_\nu X^\nu\in
E[X_1,\ldots,X_n]$ ($\nu=(\nu_1,\ldots,\nu_n)\in\nat^n$,
$X^\nu=X_1^{\nu_1}\cdots X_n^{\nu^n}$), we have 
$$v(f(\bar a))=\min_\nu \{v(b_\nu)+v({\bar a}^\nu)\}.$$In particular, $a_1,\ldots,a_n$ are algebraically independent over
$E$. Quantifier elimination then implies that
$tp(\bar a/E)$ is completely determined by
$tp(v(a_1),\ldots,v(a_n))$ (in $\Gamma$), i.e., by the set of formulas
$\sum_{i=1}^n m_i\xi_i>0$ (where $m_i\in\zee$) satisfied by
$v(\bar a)=(v(a_1),\ldots,v(a_n))$.\\[0.05in]
Let $\alpha\in E(\bar a)^s$. Then
$v(\alpha)$ belongs to the $\rat$-vector space generated by
$v(a_1),\ldots,v(a_n)$; thus there are $\rat$-linear combinations
$t_1(\bar\xi), \ldots, t_k(\bar \xi)$ such that
the values of the conjugates of $\alpha$ over $E(\bar a)$ are in
the set $\{t_1(v(\bar a)),\ldots,t_k(v(\bar a))\}$. Let $P(\bar a,X)$ be
the minimal polynomial of $\alpha$ over $E(\bar a)$.  Since
$tp(v(\bar a))\vdash tp(\bar a/E)$, there is a finite conjunction $\psi(\bar \xi)$ of
formulas of the form $\sum_{i=1}^n m_i\xi_i>0$ such that if $\bar b\in K$ is
such that $v(\bar b)$ satisfies $\psi$, and if $\beta$ is a root of $P(\bar
b,X)$, then $v(\beta)\in \{t_1(v(\bar b)),\ldots,t_k(v(\bar b))\}$.

\para\vlabel{val4}{\bf Algebraically closed valued fields with
  value group $\ree$}. Let
$E$ be an
algebraically closed field of characteristic $p>0$, and consider the
(algebraically closed) valued field $(E((t^\ree)),v)$ with its natural
valuation (see \ref{val1}). We fix a positive
integer $n$ and define $\calr$ to
be the set of $n$-tuples of $\ree^n$ which are $\rat$-linearly
independent.\\[0.05in] 
Let $\gamma=(\gamma_1,\ldots,\gamma_n)\in \calr$, and let
$\Gamma=\langle \gamma\rangle$ be
the subgroup of $\ree$ generated by the elements of $\gamma$. Let $\bar
a=(t^{\gamma_1},\ldots,t^{\gamma_n})$, let  $\alpha\in
E(\bar a)^s$,
 $P(\bar a,X)$ its minimal polynomial over $E(\bar a)$. An {\em open
   ball containing} $\gamma$ is a set $B(\gamma;\varepsilon)=\{\delta\in\ree^n\mid
\|\delta-\gamma\|<\varepsilon\}$ for some $\varepsilon>0$ (here $\|\ \|$
denotes the usual norm in Euclidean space). We will show the
 following:%\\[0.05in]
\begin{lem}\vlabel{val6} Assumptions and notation as above. 
 If $\alpha\in E((t^\Gamma))$, then there is an open ball $B$ containing
   $\gamma$ such that if $\bar b$ is an $n$-tuple in $E((t^\ree))$ with
   $v(\bar b)\in B$, then the polynomial $P(\bar b,X)$ has a root
   $\beta$ which generates an immediate extension of $E(\bar b)$. \\
Furthermore, if $Q(\bar
      X,\bar Y)\in E[\bar X,\bar Y]$ is such that $\alpha\rest 0=Q(\bar a,{\bar a}\inv)$,
      then $B$ can be chosen so that in addition the
     above  root $\beta$  satisfies $\beta\rest 0=Q(\bar b,{\bar
        b}\inv)$.

\end{lem}

% \item {\em If $\alpha\notin E((t^\Gamma))$, then there is an $\call_{\rm
%       val}(E)$-formula $\theta_{\rm ram}(\bar x)$ satisfied by $\bar a$
%     and such that 
% \end{itemize}
%
\prf The main difficulty is to say in a first order way that the
extension is immediate. We will use repeatedly Proposition
\ref{val2}. Every element of $\Gamma$ is a  $\zee$-linear
combination of elements of 
$v(\bar a)$. Hence, there is $R[\bar X,\bar Y]\in
E[\bar X,\bar Y]$ such that if $c=R(\bar a,{\bar a}\inv)$, then
$v(c)\geq 0$, and the
minimal (monic) polynomial $F(\bar a,T)$ of $\alpha'=c\alpha$ over
$E(\bar a)$ has its 
coefficients in $E[t^{\Gamma_{\geq 0}}]$. % (Note that $F(\bar a,T)=P(\bar
% a,cT)/c^{\deg(P)}$).
By Proposition \ref{val2}, we know that $$\alpha'=\lim_{m\rightarrow
  \infty}\alpha'\rest m,$$ and therefore $$\lim_{m\rightarrow
  \infty}v(F(\bar a,\alpha'\rest m))=+\infty.$$ As $F'(\bar
a,\alpha')\neq 0$ and $v(F'(\bar a,\alpha'))\geq 0$, for all sufficiently large 
$m\in\nat$ we will have $$v(F'(\bar a,\alpha'\rest m))=v(F'(\bar
a,\alpha'))\geq 0 \hbox{ and } v(F(\bar a,\alpha'\rest m))>2v(F'(\bar
a,\alpha'\rest m))+v(c).$$
We choose such an  $m\geq v(c)$, so that $v(\alpha-
\alpha'\rest mc\inv)\geq m-v(c) \geq 0$. Let $Q_m[\bar X,\bar Y]$ be such
that $\alpha'\rest m=Q_m(\bar a,{\bar a}\inv)$, and 
consider the  $\call_{\rm val}(E)$-formula $\theta_{\rm imm}(\bar x)$ which
expresses the following: \\
(i) $v(F(\bar x,Q_m(\bar x,{\bar x}\inv)))>2v(F'(\bar x,Q_m(\bar x,{\bar
  x}\inv)))+v(R(\bar x,{\bar x}\inv))$.\\
(ii) $v(Q_m(\bar x,{\bar x}\inv)R(\bar x,{\bar x}\inv)\inv - Q(\bar
x,{\bar x}\inv))\geq 0$.\\
(iii) All monomials occuring in $Q(\bar x,{\bar x}\inv)$ have valuation
$<0$, and  $v(F'(\bar x,Q_m(\bar x,{\bar
  x}\inv)))\geq 0$, $v(R(\bar x,{\bar x}\inv))\geq 0$. \\[0.05in]
Then $\bar a$ satisfies this formula. Recall that because
$\gamma\in\calr$,  $tp(\bar a/E)$ is
axiomatised by formulas of the form $\sum _i m_i v(x_i)>0$, where the
$m_i$'s are in $\zee$. Hence, there is a definable open subset $X$ of
$\ree^n$ containing $\gamma$, and such that if $v(\bar b)\in X$, then
$\bar b$ satisfies $\theta_{\rm imm}$. 
Because
$\gamma\in\calr$, there is an open ball $B$ containing $\gamma$ and
contained in this set $X$.

Let $\bar b$ satisfy $\theta_{\rm imm}$. Then, by Hensel's Lemma, there
is a root $\beta'$ of $F(\bar b,T)=0$ in the Henselization  of $E(\bar
b)$, and which furthermore satisfies
$$v(\beta'- Q_m(\bar b,{\bar b}\inv))=v(F(\bar b,Q_m(\bar b,{\bar
  b}\inv)))-v(F'(\bar b,Q_m(\bar b,{\bar b}\inv)))\geq v(R(\bar b,{\bar
  b}\inv))\geq 0$$ and $$v(\beta'R(\bar
b,{\bar b}\inv)\inv -Q(\bar b,{\bar b}\inv))\geq 0.$$
So, if $\beta=R(\bar b,{\bar b}\inv)\inv\beta'$, then $P(\bar b,\beta)=0$, $E(\bar
b,\beta)$ is an immediate extension of $E(\bar b)$, and $v(\beta-Q(\bar
b,{\bar b}\inv))\geq 0$ (i.e.: $\beta\rest 0=Q(\bar
b,{\bar b}\inv)$).

\begin{lem}\vlabel{val5} Let $A\in {\rm GL}_n(\rat)$, with $n\geq 1$. Suppose
that for every
$m\geq 1$, the characteristic polynomial of $A^m$ is irreducible
over $\rat$. Let $S\subset \calr$ with the following properties:

--- $A(S)\subseteq S\neq\emptyset$,

--- if $r\in\ree^{>0}$, then $rS=S$.

--- $S$ is closed in $\calr$.

\noindent
Then there is $\gamma\in S$ such that for every
$\varepsilon>0$, there are infinitely many integers $m$ such that
$$\Bigl\|{A^m(\gamma)\over \|A^m(\gamma)\|}-{\gamma\over
\|\gamma\|}\Bigr\|<\varepsilon.$$
\end{lem}

\prf  Our assumption on $A$ implies that for every $m\geq 1$, all
eigenvalues of $A^m$ (in $\rat^{alg}$) are distinct, and that 
$\ree^n$
has no $A^m$-stable
subspace defined over $\rat$ (other than $(0)$, $\ree^n$). Hence,  if $V$ is a proper subspace of $\ree^n$ defined over $\rat$,
then $V$ cannot contain a (non-empty) subset $D$ other than $\{0\}$ which is $A^m$-stable for
some $m\geq 1$: otherwise $\bigcap_nA^{mn}(V)$ would be a proper
$A^m$-stable subspace of $\ree^n$ containing $D$ and defined over
$\rat$.\\ [0.05in]
We will work in  the sphere
$\Ss^{n-1}=(\ree^n\setminus\{0\})/\ree^{>0}$, which we identify with 
a compact subset of $\ree^n$. We denote by $\tilde
S$, $\tilde \calr$,
the images of $S$ and $\calr$ in $\Ss^{n-1}$, i.e. with the above
identification,  $\tilde S=S\cap \Ss^{n-1}$, $\tilde \calr=\calr\cap\Ss^{n-1}$.  Let $T$ be the closure of
$\tilde S$ in $\Ss^{n-1}$; then $T$ is compact and
$A$-invariant. We call a point $c$ of $\Ss^{n-1}$ {\it recurrent}
if for every open set $U$ containing $c$ the set $\{m\in\nat\mid 
A^m(c)\in U\}$ is infinite. So we want to find a recurrent point which
is in $\tilde S$. \\
Since $S$ is closed in $\calr$, also $\tilde S$ is closed in $\tilde
\calr$, and it suffices to show 
that $T$ contains a recurrent point which is in $\tilde\calr$, 
since this point will
necessarily be in $\tilde S$.\\  {We recall the usual proof from topological dynamics:}
Take a maximal chain of non-empty closed $A$-invariant subsets of $T$;
then the intersection $C$ of this chain is non-empty (by compactness of
$T$) and is minimal closed
$A$-invariant. Thus if $c\in C$ and $k\geq 1$, then $c$ belongs to the
closure of $\{A^m(c)\mid m\geq k\}$, so that $c$ is recurrent. Note
also that $A(C)=C$.\\
Assume that $C\cap \tilde\calr $ is empty. Thus $C$ is covered by a
union of hyperplanes of $\ree^n$ which are defined over
$\rat$. There are countably many such hyperplanes, and by Baire's lemma,
there is a hyperplane $H$ defined over $\rat$ such that $C\cap H$ has
non-empty interior in $C$ for the induced topology. I.e., there is an
open subset $O$ of $\Ss^{n-1}$ such that $\emptyset\neq O\cap C\subseteq H$. 
% By Baire's lemma, there is an open subset $O$ of $\proj^{n-1}(\ree)$
% and a hyperplane $H$ defined over $\rat$ such that $C\cap O\subset
% H$. 
Note that $\bigcup_{m\in\zee} A^m(O)$ is an open set which is
$A$-invariant and intersects $C$; by minimality of $C$,  it contains
$C$, and by compactness of $C$, there is a finite subset $I$
of $\zee$ such that $C=\bigcup_{m\in I}A^m(C\cap O)$. Hence
$C\subset \bigcup_{m\in I}A^m(H)$. From $A(C)=C$, we deduce that
$\bigcap_{n\in\nat}A^n(\bigcup_{m\in I}A^m(H))$ is non-empty and
$A$-invariant. I.e, there are subspaces $V_1,\ldots,V_r$ of
$\ree^n$ which are defined
over $\rat$ and are 
 such that $C\subset V_1\cup\cdots\cup
V_r=A(V_1\cup\cdots \cup V_r)$. Thus $A$ permutes the spaces $V_i$,
and $A^{r!}(V_i)=V_i$ for any $i$. This contradicts our assumption.

\section{The main result.}

Notation and conventions are as before. In particular, $p$ denotes
$char(\Omega)$ if it is positive, and $1$ if it is $0$.

\begin{prop}\vlabel{ram3b}  Let $A$ be an abelian variety
defined over
$E=\acl(E)$,  let $B$ be a
definable modular subgroup of $A(\Omega)$, with $\evsu(B)=1$,
let $a$ be a generic of $B$ over $E$, and let $M$ be a finite $\si$-stable
Galois
extension of $E(a)_\si$. Assume that for every $\ell>1$,
$\trdeg(E(a)_\si/E)=\trdeg(E(a)_{\si^\ell}/E)$. Then for some $N$ and $b\in A(\Omega)$ satisfying $[N]b=a$,
$M\subset E(b)_\si$.
\end{prop}

\prf $B^0$ has finite index in $B$, so at the cost of replacing $a$ by
some multiple $[m]a$ (see Lemma \ref{lem1}), we may assume that $a\in B^0$. 
Replacing $a$ by $p_n(a)=(a,\ldots,\si^{n}(a))$ for some $n$ (and $A$ by $B_{(n)}$, $B$ by $p_n(B)$, see
\ref{gp1} for the notation), we may assume
that $\si(a)\in E(a)^{alg}$, and that $M=LE(a)_\si$, where $L$ is finite
Galois over $E(a)$, linearly disjoint from $E(a)_\si$ over $E(a)$, and such
that $L(\si(a))=\si(L)(a)$ (see Lemma \ref{lem4}).  Indeed, our
assumption on the difference fields $E(a)_\si$ and $E(a)_{\si^\ell}$
implies that $\evSU(p_n(B))=1$ for all $n\geq 1$ (see Remark
\ref{lem6}). Note that  $B$ is Zariski dense in $A$ (by definition of $B_{(n)}$).  \\
Let $f:W\to A$ be the
normalisation of $A$ in $L$. \\[0.05in]
By Proposition \ref{ram0}, $f$ is unramified. 
By
a result of Lang-Serre (Theorem 2 in \cite{[LS]}), this implies that $W$
is isomorphic over $E$ to an 
abelian variety $A'$, and that $f$ is a translate of an isogeny $A'\to A$.
As $E$ is algebraically closed, this implies that
$L\subset E(b)$ where $[N]b=a$ for some $N>0$.

\para\vlabel{sketch-1}  The toric analogue, Proposition \ref{prop1},  is more involved, especially in characteristic $p>0$;
before plunging into the proof, we sketch it in a
simplified setting.  Using   additive notation, up to finite index, the subgroups in question can be written:
$$B = \{x \in \gm^n\mid  \si(x)=Mx\},$$ where
$M \in {\rm GL}_n(\rat)$, and denominators
are cleared in the obvious way.  We treat here the case
$M \in {\rm GL}_n(\zee)$; a slight extension will work for $M \in {\rm GL}_n(\zee_p)
\cap {\rm GL}_n(\rat)$, but the general case is   harder.  We also consider, for simplicity,
only Galois extensions of order $p$, where $p$ is the characteristic.

We work over an
algebraically closed difference field $E$, and define
$E(B)$, the function field of $B$, to be $E(a)$ with the given
action of $\si$: $a\mapsto Ma$, where $a=(a_1,\ldots,a_n)$ is an
$n$-tuple of independent elements.    Let $L$ be a  difference field extension of $E(B)$,
which is also a Galois extension  of order $p$.  We will show that if any such $L$
exists, then $B$ is non-orthogonal to the fixed field; indeed for some $k \geq 1$, $B$ has a subgroup 
isogenous to a subgroup of 
 $\gm(\Fix(\si^k))$.  % The field $E(B)$ can be viewed as the function
 % field of the algebraic variety $(\proj^1)^n$. As such, Proposition
 % \ref{ram0}, together with the discussion on ramification,  imply that
 % the divisorial valuations on $E(B)$ which ramify in $L$ correspond to
 % subvarieties 
 % of the form $x_i=0$, or $x_i=+\infty$; furthermore, as $(\proj^1)^n$
 % does not have any proper finite unramified cover, the set of these
 % valuations is non-empty. At the cost of replacing $a$ by $a\inv$, we
 % may assume that at least one of them is ramified at $x_1=0$. If $L$

Any group embedding $\rho: \zee^n \to \ree$ % = \Hom(\gm^n,\gm)$
  induces a
Berkovich point, i.e. an $\ree$-valued valuation on $E(B)$ over $E$; namely, the   valuation satisfying
 $$v(\sum_{\nu \in \zee^n} e_\nu a^\nu   )= \min_{ e_\nu \neq 0}  \rho(\nu).$$  
Let $X$ be the space of all such valuations; so $X\subset
\Hom(\Zz^n,\ree) =\ree^n$, and  $X$ can be identified with 
$\calr$ (see \ref{val4}).     The automorphism $\si$
of $E(B)$ induces an automorphism $\si$ of $X$.   We let $S$ be the set
of valuations in $X$ which ramify in $L$.\\[0.05in] 
For any valuation $v$ in $S$, if $L=E(B)(\alpha)$ with $\alpha^p-\alpha = c \in
E(B)$ such that $v(c)<0$ and $v(c)$ is not divisible by $p$ in the value group of $v$, then it is easy
to see that $v(c)$ is uniquely determined in terms of $L,v$: any $c'$ with the same
conditions will have $v(c')=v(c)$.   (In fact, $c(1+\calo_v)$ is similarly
determined.) See \ref{ram-2} and \ref{AS}. Thus we
may define $\theta(L,v)=v(c)$. 
Moreover, any Artin-Schreier extension of $E(B)$ in which $v$ ramifies
can be described in this fashion. \\
Furthermore, the variety
$(\proj^1)^n$ has no unramified covers, and using Lemma \ref{ram0} and
the inclusion $\gm^n\subset (\proj^1)^n$, it
follows that  $S$ is
non-empty. It also satisfies the hypotheses of Lemma \ref{val5}. \\[0.1in]
We now work in the positive projectivization
 $\tilde S = (S / \ree^{>0})$ of $S$ and find a  point $\bar{v}$ such that
$\si^k(\bar{v})$ comes arbitrarily close to $\bar{v}$ (in $\tilde
S\subset \Ss^{n-1}$).  Say $\bar{v}$  is represented by $v \in S$. \\  
 Let $\Gamma$ be the
value group of $v$, a subgroup of $\ree$ isomorphic to $\zee^n$.  Then $\si^m(v) = v \circ
\si^{-m}$ is a valuation with value group $\Gamma$; we have a group automorphism
of $\Gamma$, still denoted $\si$,  satisfying
$v(\si(x)) = \si(v(x))$ (and with associated matrix $M$ with respect to
the basis $\{v(a_1),v(a_2),\ldots,v(a_n)\}$ of the $\zee$-module $\Gamma$).  \\
Of course, $\si$ does not preserve the
ordering on $\Gamma$.  But it does preserve (non)divisibility by $p$.
Moreover
since $v$ was chosen recurrent, for any given element $b$ of $E(B)$, if $v(b)<0$ then  for
infinitely many $k$, taking $\si^k(v)(b)$ close enough to $v(b)$, we have also $\si^k(v)(b)<0$.
It follows that $\theta(L,v) = \theta(L,\si^k(v))$ for infinitely many
$k$.\\
However, this means that   $\si^k$ has a  fixed point in its action on $\Gamma$, namely $\theta(L,v)$.  
Thus the characteristic polynomial of $M$ has a root $\omega$ with $\omega^k=1$.    It follows that $B$ is (up to isogeny)
a direct sum of subgroups, one of which has the form $\si(x)=M'x$ with $M'$ having a cyclotomic characteristic polynomial.  
   This finishes the sketch of
proof in this easy case.

\begin{prop} \vlabel{prop1}  Let $B$ be a definable modular
subgroup of
$\gm(\Omega)$ of $\evSU$-rank $1$, defined over  $E=\acl(E)$, and let $a$ be a generic of
$B$ over $E$. Assume that for every
$\ell\geq 1$, $\trdeg(E(a)_{\si^\ell}/E)=n$. % If $char(\Omega)=p>0$, we assume in addition that if  $n=\trdeg(E(a)_\si/E)$, then
% $\si^n(a)\in E(a,\si(a),\ldots,\si^{n-1}(a))^s$, and if $n=1$, thant
% $a\in E(\si(a))^s$.
If $M$ is a finite  $\si$-stable Galois
extension of $E(a)_\si$, then $M\subseteq E(a^{1/N})_\si$ for some
$N$.\end{prop}

\prf Our assumption on the transcendence degrees implies that
$\SU(a/E)[\ell]=1$ for every $\ell\geq 1$. By Lemma \ref{lem5}, we may 
replace $B$ by its $\si^\ell$-closure if
necessary, and 
assume that $B$ (or $B(\ell)$) is connected. Then $a$ satisfies an equation
$\prod_{i=0}^n \si^i(x^{d_i})=1$, where the
$d_i$'s are integers. We take such an equation of minimal complexity,
i.e., $n$ is minimal, and the $d_i$'s have no common
divisor (this latter condition is possible since $B$ is connected; it
will only be used in step 1, not in the rest of the proof). The
minimality of $n$ implies that $p_n(B)$ is Zariski dense in
$\gm^n$. Moreover, the polynomial $f(T)=\sum d_iT^i$ is irreducible (in
$\rat[T]$, see
\ref{end-3}). \\[0.05in] 
Since the existence of such an $M$ only depends on $qftp(a/E)$,
we will also assume that for every $m>1$, $a$ has an $m$-th root in
$B$. \\[0.05in]
{\bf Step 1}. We may assume that $\si^n(a)\in
E(a,\ldots,\si^{n-1}(a))^s$, and if $n=1$, that $a\in E(\si(a))^s$. \\ 
If $n=1$, then $p$ may divide one of $d_0$, $d_1$ but not both. Let
$m=v_p(d_1/d_0)$ ($v_p$ the $p$-adic valuation on $\rat$), and consider
$\tau=\frob^m\si$. If $m>0$, then $a$ satisfies
$a^{d_0}\tau(a)^{p^{-m}d_1} =1$, and if $m<0$, $a$ satisfies
$a^{d_0p^m}\tau(a)^{d_1}=1$. If $m=0$, then both $d_0$ and $d_1$ are
prime to $p$ and we do nothing. \\
Assume now that $n>1$, that $p$ divides $d_n$, and let $m$ be minimal
such that $\ell:=v_p(d_n)-nm=\inf\{v_p(d_i)-im\mid i=0,\ldots,n\}$. If
$\tau=\frob^m\si$, then $\tau(a)$ is in $\Ker s(\tau)$, where
$s(T)=p^{-\ell}\sum_{i=0}^n p^{-im}d_iT^i\in\zee[T]$, is irreducible (see Lemma
\ref{end-4}, and Remark \ref{lem6}), and has leading coefficient not
divisible by $p$. \\
By Lemma \ref{lem5} (see the proof), $M$ gives rise to a finite  
$\tau$-stable Galois extension of $E(a)_\tau$, so we may replace $\si$ by
$\tau$.  \\[0.05in] 
{\bf Step 2}. We may assume that for any $N>1$, the  extensions $M$ and $E(a^{1/N})_\si$ are
linearly disjoint over $E(a)_\si$.\\
It suffices to take $m$ sufficiently large, and replace $a$ by
$a^{1/m}$. \\[0.1in]
We want to show that $M$ is not
proper, i.e., that $M=E(a)_\si$. We assume this is not the case, and  define an action of $\si$ on $\rat^n$ and on $\ree^n$ by
$\si(\gamma)=A\gamma$
where 

%  $$\left( \begin{array}{cccc}
%  0 & 0 & \dots& -d_0/d_n\\
%  0 & 1 & \cdots&-d_2/d_n \\
%  \vdots&\vdots&\ddots&\vdots\\
%  0&0&\cdots&-d_{n-1}/d_n\end{array} \right)$$

$$A=\left( \begin{array}{ccccc}
0&1&0&\dots &0\\
0&0&1&\dots&0\\
\vdots&\vdots&\vdots&\ddots&\vdots\\
0&0&0&\ldots&1\\
-d_0/d_n&-d_1/d_n&-d_2/d_n&\dots& -d_{n-1}/d_n
\end{array}\right)$$
% $$A=\pmatrix{0&0&\dots&-d_0/d_n\cr
% 1&0&\cdots &-d_1/d_n\cr
% 0&1&\cdots&-d_2/d_n\cr
% \vdots&\vdots&\ddots&\vdots\cr
% 0&0&\cdots&-d_{n-1}/d_n}.$$
\noindent
Then $f(T)=\sum_{i=0}^n d_iT^i$ is the characteristic polynomial of
$A$. 
Since $\evSU(B)=1$, $f(T)$ is irreducible over $\rat$ (see
\ref{end-3}). Also,  the 
characteristic polynomial of $A^m$ for any $m\geq 1$ is irreducible over
$\rat$: otherwise, the $\si^m$-closure $B(m)$ of $B$ would contain some
definable infinite subgroup of infinite index, and therefore we would
have $\SU(B(m))>1$. Then $\evSU(B)=\evSU(a/E)=1$ would imply that
$\trdeg(E(a)_{\si^m}/E)<n$, which contradicts our assumption. 
Hence $\rat^n$ has no
proper $\si^m$-stable subspace for any $m\geq 1$.
Moreover, by modularity of $B$,
$f(T)$ is relatively prime to all polynomials of the form $T^m-p^\ell$
where $m\geq 1$, $\ell\in\zee$ (\ref{end-3}).\\[0.05in] 
Observe also that if $v$ is a valuation on $E(a)_\si$ which is trivial
on $E$ and such that
$\gamma=(v(a),\ldots,v(\si^{n-1}(a)))^T\in\ree^n$,  and if $\ell\in\zee$, then
$A^\ell\gamma=
(v(\si^\ell(a)),\ldots,v(\si^{\ell+n-1}(a)))^T$. \\[0.1in] 
{\bf Step 3}. We may assume that conjugation by $\si$ induces the identity on
$\gal(M/E(a)_\si)$.  \\
There is some $\ell$ such that $\si^\ell$ commutes with
the elements of $\gal(M/E(a)_\si)$. The divisible hull (in $\gm(\Omega)$)
of the multiplicative subgroup generated by the elements $\si^i(a)$,
$i\in\zee$, is isomorphic to $\rat^n$, and $\si$ acts on it in the
obvious manner, sending $\si^i(a)$ to $\si^{i+1}(a)$. Since $\rat^n$ has no proper
$\si^\ell$-stable $\rat$-subspace, this implies that
$\si(a),\ldots,\si^{n-1}(a)$ belong to the divisible hull of the
subgroup of $G(\Omega)$  generated by $a,\si^\ell(a),\ldots,\si^{\ell(n-1)}(a)$;
hence, for some $N\geq 1$ they belong to the group generated by
$a^{1/N},\si^\ell(a^{1/N}),\ldots,\si^{\ell(n-1)}(a^{1/N})$. We replace $a$
by $a^{1/N}$, and $\si$ by $\si^\ell$. If $char(\Omega)=p>0$, we need to
check that the separability assumptions made in Step 1 are still verified: it is
obvious when $n=1$, and follows from 
Lemma \ref{end-4}(3) when $n>1$.\\[0.05in]
%\smallskip\noindent
%\\
% The second assertion is immediate: fields between $M$ and $E(a)_\si$
% which are stable under $\si$ correspond to subgroups of
% $\gal(L/E(a)_\si)$  which are stable under conjugation by $\si$. 
%
{\bf Step 4}. We may assume that $M=LE(a)_\si$, where $L$ is finite
Galois over $E(a,\ldots,\si^{n-1}(a))$, linearly disjoint from
$E(a)_\si$ over $E(a,\ldots,\si^{n-1}(a))$,
and satisfies $L(\si^n(a))=\si(L)(a)$.\\
By Lemma \ref{lem4}, there is some $m\geq n$ such that the desired conclusion 
holds with $m$
replacing $n$. Note that $\si^n(a)\in
E(a^{1/d_n},\ldots,\si^{n-1}(a)^{1/d_n})$.
Let $N=d_n^{m-n}$: the assumption that $M\cap
E(a^{1/N})_\si=E(a)_\si$ implies that by replacing $a$ by $a^{1/N}$ and
$L$ by $L(a^{1/N},\ldots,\si^{n-1}(a^{1/N}))$,  we
obtain the  desired conclusion. \\[0.05in]
{\bf Step 5}. If $e$ is the prime-to-$p$ divisor of
$[L:E(a,\ldots,\si^{n-1}(a))]$, then  we may assume that $L$ is defined over
$E(a^e,\ldots,\si^{n-1}(a^e))$, i.e., that $L=L'(a,\ldots,\si^{n-1}(a))$
for some finite  Galois extension $L'$ of
$E(a^e,\ldots,\si^{n-1}(a^e))$, with $L'E(a^e)_\si$ finite $\si$-stable
over $E(a^e)_\si$, and $L'$ linearly disjoint from $E(a)_\si$ over $E(a^e)$.\\
We just replace $a$ by $a^{1/e}$, and $L$ by
$L(a^{1/e},\ldots,\si^{n-1}(a^{1/e}))$.\\[0.05in]
{\bf Step 6}. Consider the set $\calv$ of $E$-valuations $v$ on
$E(a,\ldots,\si^{n-1}(a))$ such that
$v(a)=\gamma_0,\ldots,\\v(\si^{n-1}(a))=\gamma_{n-1}$ are $\rat$-linearly
independent, and $v$ ramifies in $L$. Then $\calv$ is non-empty.\\
%
% For $i=0,\ldots,n-1$, consider
% $F_i=E(a,\ldots,\widehat{\si^i(a)},\ldots,\si^{n-1}(a))^{alg}$. Then
% $\bigcap_i F_i(\si^i(a))=E(a,\ldots,\si^{n-1}(a))$, so that for some $i$, we have
% $L\not\subseteq F_i(\si^i(a))$. 
We identify the tuple $(a,\ldots,\si^{n-1}(a))$ with a generic point
$(x_0,\ldots,x_{n-1})$ of
$(\proj^1)^n$. It is known that $\proj^1$ has no {\em unramified cover}, i.e,
that if $f:V\to \proj^1$ is finite, then $f$ is ramified. Hence, if $W$ is
the normalisation of $(\proj^1)^n$ in $L$, then the associated map $f:W\to
(\proj^1)^n$ is ramified.   Let $\cals$ be its ramification locus, $S$ an
irreducible component of $\cals$. Then by
\ref{ram0} and \ref{Zar}, $S$ does not intersect $\gm^n$ and is of
codimension $1$. It follows
that $S$ is of the form $x_i=0$ or $x_i=\infty$ for some $i$. Fix such
$S$ and $i$, let $v_S$ be the associated valuation.  
%
% Let $i$ be such that if
% $F=E(a,\ldots,\si^{i-1}(a),\si^{i+1}(a),\ldots,\si^{n-1}(a))$ then
% $L\not\subset F^{alg}(\si^i(a))$.
% Because the projective line has no
% proper \'etale cover, there is a valuation $w$
% on $F_i(\si^i(a))$, which is trivial on $F$,  with
% value group 
% isomorphic to $\zee$ and which ramifies in $F_iL$. Then the restriction
% of $w$ to $E(a,\ldots,\si^{n-1}(a))$ ramifies in $L$ (by \ref{ram-1}). 
% By \ref{ram0} (applied to $\si\inv$), we know that $w(\si^i(a))=\pm
% 1$: otherwise there is some irreducible polynomial $P(T)\in F_i[T]$ not
% equal to $T$, such that
% $w(P(\si^i(a)))=1$; this polynomial then yields a hypersurface of
% $\gm^n$ over which $f$ ramifies,
% and which is not contained in $\gm^n$,
%contradicting \ref{ram0}.
Let $\Delta$ be an ordered abelian group  generated by elements
$\gamma_0,\ldots,\gamma_{n-1}$ which are
$\rat$-linearly independent, and consider $\zee\oplus \Delta$ with
the lexicographical ordering. Define a valuation $v$ on
$E(a,\si(a),\ldots,\si^{n-1}(a))$ by setting $v(\si^j(a))=(0,\gamma_j)$ if
$j\neq i$, and $v(\si^i(a))=(v_S(\si^i(a)),0)$. Because the values of
$a,\si(a),\ldots,\si^{n-1}(a)$ are $\rat$-linearly independent, this
defines $v$ uniquely, and  $v$ ramifies in
$L$. 

\medskip\noindent
{\bf Step 7}. If $char(E)=0$, then $M=E(a)_\si$. \\
Because the characteristic is $0$, all valuations are defectless. We use the notation of step 5, and fix some $v\in \calv$. Then $L'$ is
defined over $E(a^e,\ldots,\si^{n-1}(a^e))$, is linearly disjoint from
$E(a,\ldots,\si^{n-1}(a))$ over $K$, and
$L'E(a,\ldots,\si^{n-1}(a))=L$. 
The number $e$ of Step 5 equals
$[L:E(a,\ldots,\si^{n-1}(a))]$, and  $L'$ is defined over
$E(a^e,\ldots,\si^{n-1}(a^e))$. So $e(L'/E(a^e,\ldots,\si^{n-1}(a^e)))$
divides $e$; moreover, 
$e\Gamma(E(a,\ldots,\si^{n-1}(a)))=
\Gamma(E(a^e,\ldots,\si^{n-1}(a^e)))$. 
So if $v\in \calv$, then $v$ does not
ramify in  $L$ by Lemma \ref{ram-2}(6) (applied to
$(K,L,M)=(E(a^e,\ldots,\si^{n-1}(a^e)), L', E(a,\ldots,\si^{n-1}(a)))$), and this contradicts step 5,
unless $L=E(a,\si(a),\ldots,\si^{n-1}(a))$ and $M=E(a)_\si$. \\[0.1in]
{\bf For the rest of the proof, we assume that the characteristic is
  $p>0$.}\\[0.1in] 
{\bf Step 8}. The set of valuations $v\in\calv$ with value group
contained in $\ree$ is non-empty.\\
Let us write $\bar a=(a,\ldots,\si^{n-1}(a))$. Extend the valuation $v$
of Step 6 to a valuation on some algebraically 
closed field $K$ containing $E(\bar a)$, and let $u\in L$ be such that $v(u)=\sum_{j=0}^{n-1}m_jv(\si^j(a))\notin
\Gamma(E(\bar a))$ (the $m_j$'s are in $\rat$); let  $P(\bar a,T)\in E[\bar
a,T]$ the minimal polynomial of $u$ over
$E(\bar a)$. Then, in the valued field $(K,v)$, $tp(\bar
a/E)\vdash \exists y\, P(\bar x,y)=0\land v(y)=\sum_j m_jv(x_j)$. \\[0.05in] 
By elimination of
quantifiers of the theory of algebraically closed fields,  and because the elements of $v(\bar a)$ are $\rat$-linearly
independent (see the discussion in \ref{val3}), there is a formula
$\psi(\bar\xi)$
satisfied by the $n$-tuple $v(\bar a)$,
which is a conjunction of formulas of the form $\sum_{j=0}^{n-1}
\ell_j\xi_j>0$ ($\ell_j\in\zee$), and   such that whenever $\bar b=(b_1,\ldots,b_{n})$ is any
$n$-tuple in  $(K,v)$,
such that   $v(\bar
b)$ satisfies $\psi$ and belongs to $\calr$ (the set of $n$-tuples of real
numbers which are $\rat$-linearly independent), then some root $\beta$  of $P(\bar b,T)=0$ has valuation
$\sum_j m_jv(b_j)$, so that $v(\beta)\notin \langle v(\bar
b)\rangle$. Choose some $n$-tuple
$\delta=(\delta_1,\ldots,\delta_{n})\in\calr$  satisfying
$\psi$, and define the $E$-valuation $v_\delta$ on $E(\bar a)$ by setting
$v_\delta(\si^j(a))=\delta_{j+1}$. Then some conjugate $u'$ of $u$
satisfies $v(u')=\sum_jm_j\delta_{j+1}\notin \langle\delta\rangle
=v_\delta(E(\bar a))$. So, 
$v_\delta$ belongs to $\calv$ and has value group contained in
$\ree$. \\[0.05in]
{\bf Step 9}. Definition of $S\subset \calr$.\\
Let $\gamma\in\calr$, and define $v_\gamma$ on $E(\bar a)$ as in the
previous step.  We  let $S$ be the set of $\gamma\in\calr$ such that
the valuation $v_\gamma$ on $E(\bar a)$ ramifies in $L$.  \\[0.05in]
{\bf Step 10}. If $\gamma\in\calr$, then there is an open ball $B'$
containing $\gamma$ and such that $B'\cap \calr\subset S$ if $\gamma\in S$, and
$B'\cap S=\emptyset$ if $\gamma\notin S$. \\%[0.05in] 
Observe that because the residue field of $(E(\bar
a),v_\gamma)$ is algebraically closed, if the valuation $v_\gamma$ does not
ramify in $L$, then  the extension $L/E(\bar a)$ is
immediate. \\%[0.05in] 
If
$\gamma\in S$, then the reasoning made in Step 8 gives us a ball $B'$
containing  $\gamma$ and such that if $\delta=(\delta_1,\ldots,\delta_{n})\in B'\cap\calr$ and $\bar
b=(b_1,\ldots,b_n)$ with $v(b_i)=\delta_i$ for $i=1,\ldots,n$, then some
root $\beta$ of $P(\bar b,T)=0$ satisfies $v(\beta)\notin \langle
\delta\rangle$, so that $E(\bar b,\beta)$ is a ramified extension of
$E(\bar b)$ since $\Gamma(E(\bar b))=\langle\delta\rangle$. \\%[0.05in]
If $\gamma\notin S$, then $L\subset E((t^\Gamma))$ and Lemma \ref{val6} gives us the result. \\[0.05in]
{\bf Step 11}. $\si(S)\subseteq S$. \\
%
%[0.05in] 
%
Let $\gamma\in S$, and fix an extension $v$ of $v_\gamma$ to $L(\si(\bar
a))$. Then the isomorphism $\si\inv $ sends  $(E(\si(\bar a)),v)$ to
$(E(\bar a), v_{\si(\gamma)})$, and therefore it suffices to show that
the restriction of $v$ to $\si(L)$ ramifies over $E(\si(\bar
a))$. \\
Observe that by step 5,
and because $e$ is the prime to $p$ divisor of $ [L:E(\bar a)]$,  the index of
ramification of $v_\gamma$ in $L$ is a power of $p$ (see Lemma
\ref{ram-2}(6) and the discussion  in Step 7). Moreover, since the residue field of $v_\gamma$ is
algebraically closed, we have $r(L/E(\bar a))=e(L/E(\bar a))$. \\%[0.05in]   
We also know that $L\si(\bar
a)=\si(L)(\bar a)$, and that $L$ and $E(\bar a,\si(\bar a))$ are linearly
disjoint over $E(\bar a)$, $\si(L)$ and $E(\bar a,\si(\bar a))$ are
linearly disjoint over $E(\si(\bar a))$. If
$\Gamma=\langle \gamma\rangle \otimes \zee[1/d_n]$ (viewed as a subgroup
of $\ree$), then $v(E(\bar
a,\si(\bar a)))\subset \Gamma$, but $v(L)\not\subset \Gamma$ because $p$
does not divide $d_n$. If
$\si(L)/E(\si(\bar a))$ is not ramified for the valuation $v$, then it
is immediate, with 
value group contained in $\Gamma$. As $E(\bar a,\si(\bar a))$ is a
totally ramified extension of $E(\si(\bar a))$, it follows that the
value group of $\si(L)(\bar a)$ is contained in $\Gamma$ (see Lemma \ref{ram-2}). This
contradicts $\si(L)(\bar a)=L(\si(\bar a))$.  \\[0.05in]
{\bf Step 12}. Choosing an $E$-embedding of $L(\si^i(a)\mid i\geq 0)$ into $E((t^\ree))$
endowed with its natural valuation (see \ref{val4}).\\ 
Observe that if $\gamma\in S$, and $r\in \ree^{>0}$, then $r\gamma\in
S$, because the valuations $v_\gamma$ and $v_{r\gamma}$ are
equivalent. If $n>1$,  
by steps 10 and 11, the conclusion of Lemma \ref{val5}
holds. Hence there
is some $\gamma\in S$ such that if $B'$ is any open ball containing
$\gamma$, then there are infinitely many $k\in\nat$ such that
${\|\gamma\|\over \|\si^k(\gamma)\|}\si^k(\gamma)\in B'$.  \\%[0.05in] 
Fix such a $\gamma=(\gamma_0,\ldots,\gamma_{n-1})$. For $i=0,\ldots,n-1$, we identify $\si^i(a)$ with
$t^{\gamma_i}$. By induction, using the fact that for $m>0$, the equation
$\prod_{i=0}^n \si^{i+m}(a_i^{d_i})=1$ determines $\si^{n+m}(a)$ up to
possibly multiplication by a $d_n$-th root of $1$, and using the fact
that $a$ is a generic of  the connected group $B$, it follows that we may
identify $\si^i(a)$ with $t^{\gamma_i}$ for all $i\geq 0$, with
$\gamma_i=v(\si^i(a))$. This defines an $E$-embedding of $E(\si^i(a)\mid
i\geq 0)$ into
$E((t^\ree))$, which extends to  $L$. We therefore
view $L(\si^i(a)\mid i\geq 0)$ as a valued subfield of $E((t^\ree))$,
the  value group of $E(\si^i(a)\mid i\geq 0)$ being contained in
$\Gamma=\langle \gamma\rangle \otimes \zee[1/d_n]$.  \\[0.05in] 
{\bf Step 13}. The final contradiction.\\
Recall that by step 5,  the index of ramification of
$v\rest {E(\bar a)}=v_\gamma$ in $L$ 
is a power of $p$. Moreover, because $E$ is algebraically closed and
$(E(\bar a),v_\gamma)$ is defectless (by \ref{defectless}(3)), it
follows that $(L,v)$ is obtained by taking first an immediate extension
of $E(\bar a)$ (namely, $L\cap E((t^\Gamma))$ by \ref{val2}), followed by a tower of Artin-Schreier extensions. Hence there is $c\in L$
such that  
$c^p-c=\alpha\in E((t^\Gamma))$, and 
the restriction of $v$ to $E(\bar a,c)$ ramifies over $E(\bar
a)$ (but its restriction to $E(\bar a,\alpha)$ does not). \\[0.05in] 
By
Proposition \ref{val2} applied to $\Gamma$,
$\alpha\rest{0}$ is a polynomial in $\bar a, {\bar a}\inv$, and we  may
assume (by \ref{AS}) that
 $\alpha\rest{0}=\sum_{\nu\in J} c_\nu \bar a^\nu$, where $\nu$ ranges over a finite
subset $J$ of $(\zee\setminus p\zee)^n$, and the $c_\nu$ are in
$E$. Let $P(\bar a,T)$ be the minimal polynomial of $\alpha$ over $E(\bar
a)$. We now use Lemma \ref{val6}: there is an open ball $B'$ containing
$\gamma$ such that whenever  $\delta\in
B'$ and $\bar b=(t^{\delta_1},\ldots,t^{\delta_n})$, then $P(\bar b,T)$ has a root $\beta$ generating an immediate
extension of $E(\bar b)$ and such that $\beta\rest 0=\sum_J c_\nu \bar
b^\nu$, with $v({\bar b}^\nu)<0$ for each $\nu\in J$. 
 Moreover, by Lemma
 \ref{val6}, these
 properties of $\bar b$ are implied by some $\call_{\rm val}(E)$-formula
 $\theta$
 satisfied by $\bar a$, namely, $v(\bar x)\in B'$. \\[0.05in]
%
%
% Then, by Proposition \ref{val2} again, there is some
% $Q(\bar a)\in E[\bar a,{\bar a}\inv]$ which approximates well $b$, i.e.,
% is such that 
% $$v(P(Q(\bar
% a)))=v(b-Q(\bar a))>\sup\{0,2v(P'(Q(\bar a)))\}.$$ 
% %
% Consider the following formula $\varphi(\bar x)$ (with parameters in $E$):
% $$v(P(Q(\bar x)))>\sup \{0,2v(P'(Q(\bar x)))\}\land
% v(Q(\bar x)-\sum_{\nu\in J} c_\nu \bar x^\nu)> 0.$$
% Then $(E(\bar a),v_\gamma)\models \varphi(\bar a)$. Let $B$ be an open
% ball containing $\gamma$ and such that if $\delta\in B$, then $(E(\bar
% a),v_\delta)\models \varphi(\bar a)$. 
Let $k\in\nat$ be such that
${\|\gamma\|\over \|\si^k(\gamma)\|}\si^k(\gamma)\in B'$. Then $(E(\bar a),v_{\si^k(\gamma)})\models
\theta(\bar
a)$: here we use the fact that $v_\gamma$ and $v_{r\gamma}$ are equivalent,
for any $r\in\ree^{>0}$. We saw that $\si^k$ defines an isomorphism of valued fields between
$(E(\bar a),v_{\si^k(\gamma)})$ and $(E(\si^k(\bar a)),v)$. Hence
$(E(\si^k(\bar a)),v)\models \theta^k(\si^k(\bar a))$, where
$\theta^k$ is obtained from $\theta$ by applying $\si^k$ to the
parameters from $E$.
We
let $\alpha_k$ be a root of $P^{\si^k}(\si^k(\bar a),T)$ satisfying
$v(b_k-\sum_{\nu\in J}\si^k(c_\nu\bar a^\nu))>0$, and $c_k$ a solution
of $T^p-T=\alpha_k$. \\[0.1in]
We let $I$ be the set of positive integers $k$ such that
${\|\gamma\|\over \|\si^k(\gamma)\|}\si^k(\gamma)\in B'$, an infinite
set by our choice of $\gamma$ in Step 12. %% Then, for each $k\in I$, we have that $c$ and
%% $c_k$ generate the same algebraic extension over  the henselisation of
%% $E(\si(\bar a)\mid i\in\nat) (b,b_k)\subset E((t^\Gamma))$.
We will shrink $I$ successively using Ramsey's theorem.\\[0.05in] 
Each $\alpha_k$ is a field conjugate of $\si^k(\alpha)$,
and it therefore follows  (by Ramsey's theorem) that there is an
infinite subset $I_1$ of $I$ 
 such that if $k<\ell$ are in $
I_1$ then $\si^{\ell-k}(\alpha_k)=\alpha_\ell$. By Step~3, every field
between $E(a)_\si$ and $M$ is a difference subfield. By Step 4
and Lemma~\ref{lem4} (the moreover part), we obtain that $E(\si^k(\bar
a),\ldots,\si^\ell(\bar a))(\alpha_k)=E(\si^k(\bar
a),\ldots,\si^\ell(\bar a))(\alpha_\ell)$. 
Hence $E(\si^k(\bar a),\ldots,\si^{\ell}(\bar a))(c_k)$ contains
$\si^{\ell-k}(c_k)$, which is a root of $T^p-T=\alpha_\ell$, and
therefore $E(\si^k(\bar a),\ldots,\si^{\ell}(\bar a))(c_k)=
E(\si^k(\bar a),\ldots,\si^{\ell}(\bar a))(c_\ell)$. The theory of
Artin-Schreier extensions tells us that there is a non-zero $c(k,\ell)\in\ffi_p$
and $g\in E(\si^k(\bar
a),\ldots,\si^\ell(\bar a))(\alpha_k)$ such that  $\alpha_k-c(k,\ell)\alpha_\ell=g^p-g$. % Thus, by Ramsey again, there is an
% infinite subset $I_2$ of 
% $I_1$ such that if $k<\ell\in I_2$, then $c(k,\ell)=1$. 
We will now use Lemma \ref{AS}. \\[0.05in] 
Let  $D=\{\nu\cdot\gamma\mid \nu\in J\}$ ($=\Supp(\alpha)\cap (-\infty,0)$;
we use the inner product notation $\nu\cdot\gamma$ for $\sum_{i=0}^{n-1}\nu_i\gamma_i$). 
Then $\Supp(\alpha_k)\cap (-\infty,0)=\si^k(D)$,   no element
of $\si^k(D)$ is divisible by $p$ in $\si^k(\Gamma)$, and
similarly for $\Supp(\alpha_\ell)\subset \si^\ell(\Gamma)\subseteq \si^k(\Gamma)$. Recall that the elements of $\gamma$ form a basis
of the free $\zee[1/d_n]$-module $\Gamma$, so that each element of $D$
corresponds to a unique tuple $\nu\in J$. By Lemma \ref{AS}(2), there is
a (unique) 
permutation $f_{k,\ell}$ of $J$ such that for all $\nu\in J$,
$f_{k,\ell}(\nu)\cdot \si^\ell(\gamma)/\nu\cdot \si^k(\gamma)$ is a
power of $p$.  By Ramsey's theorem, there is an infinite subset $I_2$ of $I_1$ such
that if $k<\ell$ are in $I_2$, then $f_{k,\ell}=id$.  
Fix $\nu=(\nu_0,\ldots,\nu_{n-1})\in
J$, $k<\ell$ in $I_2$. As above, using the freeness of $\Gamma$,  the  existence of
$m\in\zee$ such that 
$\nu\cdot \si^k(\gamma)=p^m\nu\cdot \si^\ell(\gamma)$ translates into
the following: if $u=(\nu_0,\ldots,\nu_{n-1})$, then $uA^k=p^m
uA^\ell$. So, $A^{k-\ell}$ has an eigenvalue which is a power of
$p$. This contradicts our assumption on $B'$ (see \ref{end-3}), and
finishes the proof.

\begin{thm}\vlabel{thm1} Let $\Omega$ be a model of ACFA, let
$G$  be a semi-abelian variety defined over
$E=\acl(E)$, and let $B$ be a definable modular subgroup of $G(\Omega)$
defined over $E$.  {Then the following statements hold.}
\begin{enumerate}

\item{   $B$
is stable and stably embedded.}
\item{Every definable subset of $B^n$ is a
Boolean combination of translates of definable subgroups of $B^n$. The
definable subgroups of $B^n$ are defined over $E$.}
\item{If
$a\in B$ and $F=\acl(F)$ contains $E$, and $L$ is a finite separable 
$\si$-stable extension of $F(a)_\si$, then $L$ is contained in
$F(b)_\si$ for some $b\in G(\Omega)$ and $N$ with $[N]b=a$.}
\end{enumerate}
\end{thm}

\prf Replacing $B$ by $p_r(B)$ if necessary (see \ref{gp1} for the definition),
we will assume that if $a\in B$, then $\si(a)\in E(a)^{alg}$. Hence, for
every $m\geq 1$, $\trdeg(E(a)_\si/E)=\trdeg(E(a)_{\si^m}/E)$. For $m\geq
1$, we let $B(m)$ be % the
% connected component of
the $\si^m$-closure of $B$. Then $B(m)$ is also
modular (see \ref{mod1}). Note also that by lemma \ref{lem5}, if $a$
satisfies (3) in some reduct $(\calu,\si^m)$, then it will satisfy it in
every reduct 
$(\calu,\si^m)$. By \ref{gp1}, the connected component $B^0$ of $B$ (for the
$\si$-topology) has finite index in $B$, say $m$, and therefore
$[m]B\subseteq B^0$. By Lemma \ref{lem1}, if $a\in B$, then
$[E(a)_\si:E([m]a)]_\si$ is finite, and therefore proving the
equivalence of (1) -- (3) for $[m]a$ will give the result. We may
therefore always assume that our definable subgroups $B$ and $B(\ell)$ are quantifier-free
definable and connected. 
% Note that if $\evSU(B)=1$, then  Proposition
% \ref{ram3b} and \ref{prop1} give (3).
Let us  start with some easy
remarks.\\[0.05in]
{\bf Step 1}: (1) implies (2).\\ 
$B$ with the
induced structure is stable and modular, whence weakly normal by
\cite{[HP]}. In particular every definable subset of $B^n$ is a Boolean
combination of translates of definable subgroups of $B^n$. \\[0.05in] % \\
{\bf Step 2}: (3) implies (1).\\
Let $a\in B$, and $E\subset F=\acl(F)\subset K=\acl(K)$, and assume that
$F(a)_\si$ and $K$ are independent over $F$. By (3), we know that
all finite separable $\si$-stable extensions of $K(a)_\si$ are contained in
$\acl(Ea)K\subseteq \acl(Fa)K$. Assume that $L$ is a finite separable $\si$-stable
extension of $\acl(Fa)K$. By Lemma \ref{lem31}, there is some $K_1$
independent from $K(a)_\si$ over $F$ and such that $L\subset
\acl(K_1a)M$ for some finite separable $\si$-stable extension $M$ of
$KK_1(a)_\si$. By (3), we get  $M\subset \acl(Fa)\acl(KK_1)$, therefore
$L\subset \acl(Fa)\acl(KK_1)\cap \acl(Ka)=\acl(Fa)K$ (Because $K_1$ is
independent from $\acl(Ka)$ over $F$; see e.g. Remark 1.9(2) in
\cite{[CH]} applied to $A=\acl(Fa)$, $B=K$ and $C=K_1$). By Lemma
\ref{bab3} (applied to $K\acl(Fa)$), $tp(a/F)\cup qftp(a/K)$ is
complete, which shows that $tp(a/F)$ is stationary.  
Hence  every type over an algebraically closed set
 which is  
realised in $B$ is stationary. The result follows by \ref{sse}. %  (See Lemma 2 in \cite{[CH]}. There is a
% typo, the set $A$
% of (3) should satisfy $A={\rm acl}^{\rm eq}(A)$)
Moreover, by Lemma
\ref{lem5}, we obtain that all $B(m)$ are stable and stably
embedded. \\[0.05in] 
%
%
%
%
%
% We may assume that $B=B^0$. Replacing $B$ by $p_r(B)$ if necessary,
% we will assume that if $b\in B$, then $\si(b)\in E(b)^{alg}$. 
% If $m\geq 1$, let $B(m)$ be the
% connected component of the $\si^m$-closure of $B$. Then $B(m)$ is also
% modular (see \ref{mod1}). The idea of the proof is simple: show (3) for
% types of evSU-rank $1$, obtain (1), and then use \ref{sse2} to get (1) in
% general, from which (2) and (3) follow. Unfortunately, the toric case
% makes things a little more
% complicated in positive characteristic. Let us first make some easy remarks.\\
%
% 
%%%%%%%%%%%%%%%%%%%%%%%%%%%%%%%%%%%%%%%%%%%%%%%%%%%%%%%
%
{\bf Step 3}: (1) for all $B(\ell)$ implies (3). \\
Let $F=\acl(F)$ contain $E$, and
let $a\in B$. Let $C$ be the 
$\si$-closed connected subgroup of $B$ such that $a$ is a generic of
the coset $a+C$ over $F$. Let $a_1$ be a generic of $a+C$ which is independent from $a$
over $F$. Then $tp(a/Fa_1)$ is completely determined by the class of
$(a-a_1)$ in $C/C^*$, where $C^*=\bigcap_{m\geq 1}[m]C$. This comes
from the fact that $tp(a/Fa_1)$ is uniquely determined by the set of
cosets of
definable subgroups $D$ of $C$ which are defined over $F(a_1)_\si$ and
which contain $(a-a_1)$; as $(a-a_1)$ is a generic of $C$ over
$F(a_1)_\si$, these subgroups $D$ must have finite index in $C$, and
therefore  contain  $[m]C$ for some positive integer $m$. For each $m>1$ choose
$b_m\in G$ such that $[m]b_m=a-a_1$, and let $F_1=\acl(Fa_1)$.
Then $tp(a/F_1)$ is completely
determined by $qftp(a,b_2,\ldots,b_m,\ldots/F_1)$. Similarly, for
every $\ell\geq 1$, $tp(a/F_1)[\ell]$ is completely determined by
$qftp(a,b_2,\ldots,b_m,\ldots/F_1)[\ell]$. That is (see \ref{bab3}), the field
$F_1(a,b_2,\ldots )_\si$ has no finite proper separable $\si$-stable
extension.\\[0.05in] 
Let $L$ be a finite separable $\si$-stable extension of $F(a)_\si$. Then
$LF_1(a,b_2,\ldots)_\si$ is a
finite $\si$-stable
extension of $F_1(a,b_2,\ldots)_\si$, so that $L\subset F_1(b_m)_\si$ for
some $m$.\\[0.05in]
Note that $F_1$ contains a root of $[m]x=a_1$. Hence
$F_1(b_m)_\si=F_1(c_m)_\si$ where $[m]c_m=a$, and $L\subset
F_1(c_m)_\si$. As $a_1$ was independent from $a$ over $F$, we have
$F_1(c_m)_\si\cap \acl(Fa)=F(c_m)_\si$, so that  $L\subset
F(c_m)_\si$.\\[0.1in] 
It therefore suffices to show (3) or to show (1) for all
$\ell$. \\[0.1in] 
%
%
% Putting all these observations together, and using \ref{sse2} for $(*)$,
% we have shown:\\ 
% {\em $(*)$ If $C\subset B$ is a definable subgroup, and $C$, and $B/C$ both
%   satisfy (3), then so does $B$. \\
% $(**)$ If  $B(m)$ satisfies (3), then so
%   does $B$.}\\[0.05in]
%
% It now remains to
% show (3).  
% The proof will proceed in several additional steps:\\ 
% Show the result under the assumption that it holds for subgroups of
% evSU-rank $1$.\\
% Special cases of evSU-rank 1 subgroups.\\
% The general case of evSU-rank 1 subgroups (of tori).\\[0.1in] 
%
{\bf Step 4}. Case where $B$ has evSU-rank $1$.\\
Assume that $B$ has evSU-rank $1$. \\
If $G$ is an abelian
variety, then Proposition \ref{ram3b} gives (3). \\

% Let $a$ be a generic of $B$ over $F$. If $G$ is an abelian
% variety and $L$ is a finite separable
% $\si$-stable extension of $F(a)_\si$, then 
% Proposition \ref{ram3b} gives (3). \\
Assume now that $G$ is not abelian, but has an abelian quotient
$A$ (via a morphism $\pi$) such that $\pi(B)$ is infinite; then 
the restriction of $\pi$ to $B$ has finite kernel, and by the previous
case, we obtain that $\pi(B)(\ell)$ is stable and stably embedded for
every $\ell>0$. As
$\Ker(\pi)\cap B$ is finite, 
Proposition
\ref{sse2} gives that 
$B(\ell)$ is stable and stably embedded for every $\ell>0$, and
therefore gives also (3) and (2). 
This show (1) and (3) in the ``abelian case''. \\
Assume now that the Zariski closure of $B$ is contained in
some toric subvariety of $G$. Recall that $\si(a)\in E(a)^{alg}$. Without loss of generality, 
$G=\gm^r$ for some $r$. Let $a_1\in \gm$ be an element
of the tuple $a$ which does not belong to $E$. Then
$\trdeg(E(a_1)_\si/E)=\trdeg(E(a)_\si/E)$ (because $\evSU(a/E)=1$), and as above, the
projection $\pi$ of $\gm^r$ on the corresponding copy of $\gm$ restricts to a
morphism  on $B$ with finite kernel. 
Applying Proposition \ref{prop1} to $a_1$ shows that if $L$ is a finite
$\si$-stable extension of $F(a_1)_\si$, then $L\subset F(a_1^{1/N})_\si$
for some integer $N$. Hence, $\pi(B)(\ell)$ is stable and stably
embedded for every $\ell>0$,
and so is $B(\ell)$. This finishes the proof when 
$\evSU(B)=1$.\\[0.05in]
{\bf Step 5}. The general case.\\
Let $a$ be a generic of $B$ over $E$, and $n={\rm
  tr.deg}(E(a)_\si/E)$. Recall that  if $a\in B$, then
 $\si(a)\in E(a)^{alg}$. 
Let $m$ be
large enough so that 
$\SU(a/E)[m]=\evSU(a/E)=k$, and replace $B$ by $B(m)$ (or by its
connected component, but we will keep the notation $B(m)$). We now work in
$(\Omega,\si^m)$. Using the modularity of $B(m)$, there is a
sequence $(0)=B_0\subset B_1\subset \cdots \subset B_k=B(m)$ of
subgroups of $B(m)$, with $[B_i:B_{i-1}]=\infty$ for $i=1,\ldots,k$:  take
a sequence $F_{k-1}\subset \cdots\subset F_1$ of difference subfields
of $\Omega[m]$ such that $SU(a/F_i)[m]=i$. Then for each $i$, $a$ is the generic
of a coset of a $\si^m$-closed subgroup $B_i$
of $G$ of $SU[m]$-rank $i$, by Theorem \ref{f1}. If $C_i=B_i/B_{i-1}$,
then 
 $\evSU(C_i)=1$ for
$i=1,\ldots, k$. We will now show that each $C_i$ ``lives'' in a
semi-abelian variety,  so that it  satisfies (3) by the previous
steps. Let $\pi$ be a morphism (of algebraic groups) from 
the Zariski closure of $B_i$ (inside $G(\Omega)$) onto
some simple semi-abelian variety $H$, such that $[\pi(B_i):\pi(B_{i-1})]$
is infinite. Let $f(\si)\in \End_\si(H)$ be such that $\pi(B_{i-1})$ is
commensurable to $\Ker(f(\si))$; multiplying $f$ by an integer, we may
assume $\pi(B_{i-1})$ is contained in $\Ker(f(\si))$; then $f(\si)\circ \pi(B_i)$ is an infinite subgroup of
$H(\Omega)$, because $\Ker (f(\si))$ has infinite index in $\pi(B_i)$, and the induced map $C_i\to H(\Omega)$ has
finite kernel. 
Hence, by
Step 4 and \ref{sse2}, each $C_i$ is stable and stably embedded, and by
\ref{sse2}, this implies that $B(m)$ is stable and stably embedded. The
same reasoning shows that all $B(m\ell)$ are stable and stably embedded; 
using Step 3
and Lemma \ref{lem5},
$B(m)$ and $B$ satisfy (3), and therefore also (1) and (2). This
finishes the proof. 

\para{\bf Some remarks about $B/[n]B$.} 
Let $A$ be a semi-abelian variety
defined over $E=\acl(E)$,
and let $B\subset A(\Omega)$ be a definable modular subgroup of
$A(\Omega)$, $n$ an integer bigger than $1$. The induced structure of
the definable group $B$ is largely determined by the finite index
subgroups $[n]B$. For instance if 
$[n]B=B$ for all $n>0$, then $B$ is strongly minimal and the induced structure reduces
to a module structure  (cf. \cite{HL}).    However this only occurs in the rare event that $B$ is torsion-free;
it will usually have infinite torsion.  One thus wants to  understand the finite imaginary sort $B/[n]B$ in
terms of a similar quotient for the ambient algebraic group $A$. In what
follows, we will assume that $B$ is a quantifier-free definable subgroup
of $A(\Omega)$, of finite SU-rank.\\
Since $\Omega$ is elementarily equivalent to an
ultraproduct of difference fields  $(\ffi_p^{alg},\Frob _q)$ (see \cite{[H2]}), we
know that
$[B:[n]B]=|B\cap A[n]|$: the group $B$ is elementarily equivalent to an
ultraproduct of finite groups, and one considers the endomorphism $x\mapsto
[n]x$. We will be able to give a precise description of $B/[n]B$ in two cases: when
$B=\Ker(f)$ for some $f\in {\rm End}_\si(A)$; and when $B$ is
connected (for the $\si$-topology). Any definable group of finite
SU-rank $B$ is isomorphic (via a map $p_N$) to a finite index subgroup of some $B'$ of the first type,
and contains a finite index subgroup $B''$ of the second type, so our
description is rather complete.

\begin{prop} Let $A$ be a semi-abelian variety defined over $E=\acl(E)$, and let
$B\subset A(\Omega)$ be a quantifier-free definable  subgroup of $A(\Omega)$ of finite SU-rank, $n>1$
an integer. \begin{enumerate}
\item Assume that $B=\Ker(f)$ for some
  $f\in \End_\si(A)$. Then $$B/[n]B\simeq A[n]/f(A[n]).$$

\item Assume that $B=\tilde B^0$,  let $N$ be such that $B=\tilde
  B_{(N)}$, and consider the semi-abelian variety $D=B_{(N)}$. Take a definable homomorphism $F:D\to D'$, where
  $D'$ is the quotient of the semi-abelian variety $D^\si$ by some
  finite subgroup $C_2$, and the map $F$ is of the form
  $f(p_N(x))-g(\si(p_N(x))$ for some algebraic homomorphisms $f,g$, such
  that $B=\{b\in A(\Omega)\mid F(p_N(b))=0\}$. Then  

$$B/[n]B\simeq D'[n]/F(D[n]).$$

\end{enumerate}

\end{prop}

\prf (1) Let us first suppose that  $B=\Ker(f)$. 
Consider $A[n]/f(A[n])$. Then $[A[n]:f(A[n])]=|\Ker(f)\cap
A[n]|=|B\cap A[n]|$. Define $$\varphi:B\to A[n]/f(A[n])$$ as follows:
if $a\in B$ take $b\in A$ such that $[n]b=a$, and define $\varphi(a)=
f(b)+f(A[n])$.
Since distinct choices of  $b$ differ by an element of $A[n]$, this map is
well-defined. One checks easily that it is a group homomorphism, and that
its kernel is precisely $[n]B$, so that it defines an isomorphism
between $B/[n]B$ and $A[n]/f(A[n])$. \\[0.05in]
(2) We will use the general description of
quantifier-free definable subgroups given in \ref{gp1}, and follow its
notation. 
Replacing $B=\tilde B^0$ by $p_N(B)$ for some
$N$ and $A$ by $B_{(N)}=D$, we may assume that $B=\tilde B_{(1)}$, and that $B$
is Zariski dense in the semi-abelian variety $A$ (this is where we use
that $B$ has no quantifier-free definable subgroup of finite index:
$B_{(N)}$ is connected).
If we define $C_1,C_2$ by $$C_1\times \{0\}=B_{(1)}\cap (A\times \{0\}),
\qquad
\{0\}\times C_2=B_{(1)}\cap (\{0\}\times A^\si),$$ the group
$B_{(1)}/C_1\times C_2$ is the graph of a (definable in ACF) group
isomorphism
$$f:A/C_1\to A^\si/C_2.$$
If $h_1:A\to A/C_1$ and $h_2:A^\si\to A^\si/C_2$ are the natural isogenies, $a\in
B$ if and only if $fh_1(a)=h_2(\si(a))$. \\
Let $A'=A^\si/C_2$, let $F:A\to A'$ be defined by
$F(x)=fh_1(x)-h_2(\si(x))$. Then $B=\Ker(F)$.
Define $\varphi:B\to A'[n]/F(A[n])$ as follows: if $a\in B$, let $b\in A$
be such that $[n]b=a$, and set $\varphi(a)=F(b)+F(A[n])$. As $A$ and $A'$
are isogenous semi-abelian varieties, we know that $|A[n]|=|A'[n]|$,
whence $[A'[n]:F(A[n])]=|\Ker(F)\cap A[n]|=|B\cap
A[n]|$, and we get an isomorphism between $A'[n]/F(A[n])$ and
$B/[n]B$.

\para{\bf Algebraic dynamics}. We call {\em algebraic dynamic} a pair
$(V,\phi)$ consisting 
of a (quasi-projective, irreducible) variety $V$, together with a
dominant rational map $\phi:V\to V$. A map between algebraic dynamics
$(V,\phi)$ and $(W,\psi)$ is a dominant rational  map $f:V\to W$ such
that $f\circ 
\phi=\psi\circ f$. For $n>0$, $\phi^{(n)}$ denotes
$\phi\circ\phi\cdots\circ\phi$ ($n$ times).\\[0.05in]
If $(V,\phi)$ is defined over the field $E$, then we put the structure
of  a
(non-inversive) difference field  on $E(V)$ by setting $\si(a)=\phi^*(a)=a\circ
\phi$. Note that $\si$ is the identity on $E$. A dominant rational map
$f:(V,\phi)\to (W,\psi)$ (defined over $E$) then corresponds to an embedding of difference fields
$f^*: E(W)\subset E(V)$. 

\begin{prop}\vlabel{dyn}  Let $G$ be
  a semi-abelian variety, defined over an algebraically closed  field $E$. Let
  $\phi:G\to G$ be a dominant 
  endomorphism, and assume that the set  $\{x\in G(\Omega)\mid
  \si(x)=\phi(x)\}$ is one-based (in some existentially closed difference field
  $\Omega$ containing $(E,id)$). Let $(W,\psi)$ be an algebraic
  dynamic, and 
  $f:(W,\psi)\to (G,\phi)$ a finite separable map, everything being defined over
   $E$. Then there is a birational map
  $g:(W,\psi)\to (H,\rho)$, and an isogeny $h:(H,\rho)\to (G,\phi)$ such
  that $h\circ g=f$, where $H$ is a semi-abelian variety, and $\rho$ a
  dominant endomorphism of $H$. 
\end{prop}

\prf We work in the difference field  $\Omega$, a model of ACFA. We
take $b\in W$, generic over $E$ and satisfying $\si(b)=\psi(b)$, and let
$a=f(b)$. Then $a\in G$, is generic over $E$ and satisfies
$\si(a)=\phi(a)$. Then $E(a)\subset E(b)$ are closed under $\si$, but
not  under $\si\inv$, unless $\phi$ and $\psi$ are
isomorphisms. Moreover, $E(b)$ is a finite separable extension of
$E(a)$. Results on the  limit degrees (equal to $1$ in both cases), and
the fact that $b$ is algebraic over $E(a)$ imply that $E(b)_\si$ is a
finite (separable) $\si$-stable extension of $E(a)_\si$ (see \cite{[C]},
Theorem 5.22.XVI). \\ 
Let $B\subset G(\Omega)$ be the subgroup  defined by the difference equation
$\si(x)=\phi(x)$. 
By Theorem \ref{thm1}, $E(b)_\si\subset E(a^{1/N})_\si$
for some integer $N>0$. This implies that for some $n\geq 0$, $N>0$,
$E(b)\subset E(\si^{-n}(a^{1/N}))$. Note that $[N]\phi^n:G\to G$ is an
isogeny, which sends $\si^{-n}(a^{1/N})$ to $a$. It factors as $h_1\circ
h_2$, where $h_1$ and $h_2$ are isogenies such that 
$h_1$ is separable, and (setting $c=h_2(\si^{-n}(a^{1/N}))$)   
$E(b)=E(c)$. Let $H=h_2(G)$, $\rho$ the
endomorphism of $H$ such that $\si(c)\in E(c)$, and $h=h_1$. 
\begin{rem} The hypotheses of Proposition \ref{dyn} are not very
  friendly to non-logicians. \\
In the notation of the proof of \ref{dyn}, assume that $B$ is not
one-based. Then, there is a semi-abelian (simple) variety $H$, and an algebraic
map $h:G\to H$, such that $h(B)=C\subseteq H(\tau)$, with
$\tau=\si^m\frob^n$ for some $m\geq 1$ and $n\in\zee$, and $h(C)$ infinite. If $\tau=\si^m$,
then $H$ and $h$ can be taken defined over $E$; if $\tau=\si^m\frob^n$
with $n\neq 0$, then $H$ can be taken to be defined over some finite
field, and $h$ over $E$ as before. \\[0.05in]
In the first case, as for  $m$ sufficiently large, $\si^m$  commutes with all elements of
$\End(G)$, it will follow that for some $m$, the map $\phi^{(m)}-id$ is
  not onto. \\[0.05in]
The second case is not as easy to describe. If $G$ is simple and equals $\gm$,  then 
$\phi^{(m)}-\frob^{-n}$ is not onto (and we get that $-n>0$). Assume now that $G$ is abelian simple,
and let $h:G\to H$ be the isogeny given above,  $h^*:H\to G$  its dual,
and $M$ the integer such that $hh^*=[M]$. Then for some $m$ and $n$, we
have that $M\phi^{(m)}-h^*\frob^nh$ is not dominant.  

The general case is harder to describe. 
\end{rem}

\para{\bf Problem}.  In positive characteristic, do there exist any
stable definable subgroups of $\ga$? \\[0.1in]
It is easy to show that in characteristic $0$ there are none since any
quantifier-free definable subgroup of $\ga$ is defined by linear
difference equations. See
e.g. Theorem 5.12 in \cite{[CH]}.

\bigskip \noindent
Current addresses:\\[0.1in]
D\'epartement de Math\'ematiques et Applications (UMR 8553)\par\noindent
Ecole Normale Sup\'erieure\\
45 rue d'Ulm\\
75230 Paris Cedex 05\par\noindent
France \par\noindent
e-mail: {\tt zchatzid@dma.ens.fr}\\[0.1in]
Institute of Mathematics \par\noindent
The Hebrew University \par\noindent
Givat Ram \par\noindent
91904 Jerusalem \par\noindent
Israel \par\noindent
e-mail: {\tt ehud@math.huji.ac.il}

\end{document}